\documentclass[12pt]{article}
\title{Geodesic connectedness of multiwarped spacetimes.}
\author{Jos\'{e} Luis Flores and Miguel S\'anchez \thanks{The authors acknowledge warmly to Prof. J.L. G\'amez and Prof. D. Arcoya for having some discussions. Research partially
supported by
DGICYT grant number PB97-0784-C03-01.} \\ Depto. Geometr\'{\i}a y
Topolog\'{\i}a, Fac. Ciencias, Univ. Granada, \\ Avda. Fuentenueva s/n, 18071-Granada, Spain. \\E-mail a.: sanchezm@goliat.ugr.es}

\newtheorem{lemma}{Lemma}
\newtheorem{remark}{Remark}
\newtheorem{theorem}{Theorem}
\newtheorem{proposition}{Proposition}

\newtheorem{convention}{Convention}

\newtheorem{definition}{Definition}

\font\ddpp=msbm10 scaled \magstep 1 
\def\R{\hbox{\ddpp R}}    
\def\L{\hbox{\ddpp L}}    
\def\N{\hbox{\ddpp N}}    

\def\be{\begin{equation}}
\def\ee{\end{equation}}

\begin{document}
\maketitle

\begin{center}
{\it Running title:} Connectedness of spacetimes
\end{center}

\begin{abstract}
A new technique for the study of geodesic connectedness in a class of Lorentzian manifolds is introduced. It is  based on arguments of Brouwer's topological degree for the solution of functional equations. It is shown to be very useful for multiwarped spacetimes, which include different types of relativistic spacetimes.
\end{abstract}

\newpage 

\section{Introduction}

As far as we know, the results about geodesic connectedness in Lorentzian manifolds (i. e. the problem as to whether each pair of their points can be joined by a geodesic) can be grouped into four types, as follows:
\begin{itemize}
\item Results on space-forms, first obtained in \cite{C-M} and compiled in the books \cite[Chapter 11]{Wo} or \cite[Chapter 9]{O}. In particular, a positive Lorentzian spaceform is geodesically connected if and only if it is not time-orientable. 
\item Results on {\it disprisoning} and {\it pseudoconvex} manifolds; both geometrical concepts were introduced by Beem and Parker \cite{B-P} and are studied in the book \cite{B-Eh-Ea}.  Lorentzian manifolds satisfying these two conditions and having no conjugate points are shown to be geodesically connected; moreover, in this case, a Lorentzian Hadamard-Cartan theorem holds, which implies that the manifold is diffeomorphic to $\R{}^n$.
\item Results obtained by means of variational methods. Roughly, geodesics are seen as critical points of the (strongly indefinite) action functional defined on curves joining two fixed points, and some techniques are developed to make sure that this functional admits critical points. This method was introduced by Benci and Fortunato \cite{B-F88}, \cite{B-F90} and since then, it has proved to be very fruitful. In fact, a series of results have shown the geodesic connectedness of many families of Lorentzian manifolds which  generalize most of the classical spacetimes used in General Relativity (see the book \cite{Ma} or the more recent references \cite{G-M}, 
\cite{A-S}, \cite{G-P}, \cite{C-S}, \cite{Pi}). 
\item Results based on a direct integration of the equation of the geodesics, in especially interesting  cases from either a mathematical or a physical point of view \cite{Sa97a}, \cite{Sa98}, \cite{C-S}.
\end{itemize}
Other techniques allow the study of the existence of causal (timelike or null) geodesics between two given points  (even though in some cases the existence of a timelike geodesic can be deduced from variational methods, as shown in some of the previous references). Among the  results for these geodesics, the classical Avez-Seifert one is especially relevant: in a globally hyperbolic spacetime, two points can be joined by a causal geodesic if and only if they can be joined by a causal curve.  The aim of this paper is to introduce a new technique for the study of geodesic connectedness in certain Lorentzian manifolds, based on arguments involving topological degree of solutions of functional equations.

We will concentrate on {\it multiwarped spacetimes}, which are product manifolds $I\times F_1 \times \cdots \times F_n,$ $I\subseteq \R{}$   endowed with a  metric $g = -dt^2 + \sum_{i=1}^n f_i^2(t) g_i, t\in I$ (see Section 2 for precise definitions). From a physical point of view, these spacetimes are interesting,   first, because they include classical examples of spacetimes: when  $n=1$ they are the {\it Generalized Robertson-Walker} (GRW) spacetimes, standard models of inflationary spacetimes \cite{Sa98}; when $n=2$, the intermediate zone of Reissner-Norsdstr\"om spacetime and the interior of Schwarzschild one appear as particular cases \cite{Sa97}. Moreover, multiwarped spacetimes may also represent relativistic spacetimes together with internal spaces attached at each point (see \cite{Sc} and references therein). 

From a mathematical point of view, there are serious problems  studying geodesic connectedness of multiwarped spacetimes by using previous techniques, even for $n=2$. In fact, the results known regarding this  are (in what follows, the {\it line} at $x\in F_1 \times \cdots \times F_n$ is the set 
$L[x]=\{ (t,x)| t\in I\}$):
\begin{enumerate}
\item For $n=1,2$, causal geodesics are completely characterized in \cite{Sa97}, showing  an Avez-Seifert type result just under the assumption of {\it weak convexity} for the fibers. By weak convexity of a  Riemannian manifold we mean geodesic connectedness by minimizing (not necessarily unique) geodesics; in the cited reference,  this hypothesis is shown to be essential and it will be assumed in this paper.
\item For $n=1$, the previous result and elementary arguments on continuity and causality show \cite[Section 3]{Sa98}: if any point $z=(t,x)\in I\times F (F\equiv F_1)$ and any line $L[x']$ can be joined by a future as well as a past directed curve then the spacetime is geodesically connected; equally, the conclusion holds if $f\equiv f_1$ satisfies 
\begin{equation}
\label{f}
\int_a^{c} f^{-1} = \int^b_{c} f^{-1} = \infty. 
\end{equation}
A remarkable example where this condition is not satisfied is de Sitter spacetime, which can be seen as a GRW spacetime with $f= \cosh$, and it is  non-geodesically connected.
\item For $n=2$,  the geodesic connectedness of spacetimes with a qualitative behavior such as the Reissner-Nordstr\"{o}m Intermediate one is proven in \cite{Gi}. This result was extended in \cite{G-M} to the case when the fibers admit more general topological boundaries. Remarkably, the technique there shows the equivalence  of the action functional and a simpler functional. Thus,  a set of conditions must be imposed to ensure the existence of critical points for this functional.
In our more general setup, we can make all the fibers equal to intervals of $\R{}$ to simplify things, and, then, it is easy to check that
any geodesic $(\tau(t),\gamma_1(t),...,\gamma_n(t))$ joining $(\tau_0,x_1,...,x_n)$ and $(\tau'_0,x'_1,...,x'_n)$ can be recovered as a critical point of the functional
\begin{equation}
\label{func}
\tau \rightarrow -\int_0^1 \tau^{'2}(t) dt + \sum_i \frac{(x_i'-x_i)^2}{ \int^1_0f_i^{-2}(\tau(t)) dt } 
\end{equation}
defined on curves $\tau:[0,1]\rightarrow I $ joining $\tau_0$ and $\tau'_0$.
\end{enumerate}

In this paper we prove,

\begin{theorem}\label{t1} A spacetime $(I\times F_{1}\times \ldots \times F_{n},g)$ with weakly convex fibers is geodesically connected if the following condition holds: any point of the spacetime can be joined with any line by means of both, a future directed and a past directed causal curve. 
\end{theorem}
This condition can be expressed easily in terms of the warping functions, being equivalent to formulae (\ref{e28}) (Proposition \ref{p1}) and it is completely natural from a mathematical point of view. 
After some preliminaries in Section 2, we prove  Theorem \ref{t1} in the
following three sections. First, in
Section 3, connection by causal geodesics is characterized. Most of the  proof is completely analogous to the case $n=2$ solved in \cite{Sa97}; so, essentially, we prove Lemma \ref{l3} only, which is a non-trivial generalization of a step for  $n=2$. In Section 4 we prove a particular case of Theorem \ref{t1} (conditions (\ref{e28}) are replaced by the stronger 
(\ref{e24})) by means of topological degree. For the general result, proven in Section 5, some additional problems appear, making  both the concept of {\it fake} geodesic and the hypothesis in Theorem \ref{t1} natural.

Furthermore, there are geodesically connected multiwarped spacetimes where Theorem \ref{t1} is not applicable, some of them of special physical interest (such as the Reissner-Nordstr\"om intermediate spacetime itself). In fact, condition (\ref{f}) (or its generalization (\ref{e28})) seems appropiate  when $I=\R{}$ or $f$ goes to zero at the extremes; nevertheless,  a strip $I\times \R{}^n, I\neq \R{}$ in Lorentz-Minkowski spacetime $\L{}^{n+1}$ does not satisfy this condition.
However, Theorem \ref{t1} does not cover all the possibilities of the technique, so we will give a more general version of this theorem in the last Section, under a hypothesis which is close to a necessary condition (Theorem \ref{t5}). Thus, all previous results are reproven or extended; in particular, geodesic connectedness of Reissner-Nordstr\"{o}m Intermediate spacetime is reproven,  and natural conditions for the existence of critical points of the functional (\ref{func}) can be obtained. Moreover, the accuracy of our technique is shown by proving the geodesic connectedness of Schwarzschild inner spacetime (Theorem \ref{t6}, Remarks \ref{r6} and \ref{r0}).
It is also worth pointing out that our results can be trivially extended to the case when the fibers $F_i$ are incomplete and has a Cauchy boundary $\partial F_i$. In fact, it is enough to wonder when the structure of  $\partial F_i$  implies that $F_i$ is weakly convex (if $F_i \cup \partial F_i$ is a differentiable manifold with boundary, $F_i$ is weakly convex  if and only  if the second fundamental form of the boundary, with respect
to the interior normal, is positive semidefinite; for more general results, see \cite{B-G-S}). This improves the previous results on this case too.

In the remainder of the present section, we give an intuitive idea of the techniques in Sections 4 to 6, for $n=2$ fibers.

Fix a point $z=(\tau_0,x_1,x_2) \in I\times F_1\times F_2$ and try to connect it with a point $z'$ in the line $L[x_1',x_2']$ ($x_1'\neq x_1, x_2'\neq x_2$) by means of the geodesic $\gamma(t) = (\tau (t), \gamma_1(t), \gamma_2(t))$. If $z'$ belongs to the future or past of $z$, the problem is solved by Section 3; so, it is necessary to study just when $z'$ belongs to a compact interval $J$ of the line, such that each one of its extremes is causally related with $z$. For  a multiwarped spacetime, the projections $\gamma_1, \gamma_2$ of the geodesic $\gamma$ are pregeodesics of the fibers, that is, up to (probably different) reparameterizations, geodesics of the fibers. It is natural to only consider  the case when these geodesics  minimize the corresponding Riemannian distance (if they are not unique,  we also assume a fixed choice has been carried out). So, $\gamma$ can be characterized by three paremeters, say, $K$, directly related to $\tau'(0)$, and $c_1, c_2 $, related to the initial speed of each pregeodesic. Moreover, $\gamma$ can be reparametrized to assume $c_1 + c_2 =1$, reducing both parameters $c_i$'s to just one $c (\equiv c_1) \in (0,1)$, and the domain of $K$ can be assumed to be a compact interval $[K^-,K^+]$ such that the geodesics with $K=K^-,K^+$ are necessarily causal for all $c$.

When $\tau(t)$ is not constant, each $\gamma_i(t)$ can be reparameterized by $\tau$, up to a set of isolated points which will be specifically taken into account. Now consider the functions $s_i(c,K) \; i=1,2$  which maps each (geodesic) $c, K$ in the length of the interval $I=(a,b)$ covered by the  parameter $\tau \in I$ when $\gamma_i(t(\tau))$ goes from $x_i$ to $x_i'$ (see the beginning of the proof of Theorem \ref{t3}). Then, the zeroes of the function $\mu\equiv 1-\frac{s_{2}}{s_{1}}$ represent geodesics joining $z$ and the line $L[x'_1,x'_2]$. Under assumption (\ref{e24}) this function satisfies:
\be
\label{lim}
\lim_{c\rightarrow 0}\mu(c,K)<0, \quad \quad  \lim_{c\rightarrow 1}\mu(c,K)>0. 
\ee
 Thus, for each $K$ there exists at least one zero of $\mu$ and, by  arguments which  naturally involves Brouwer's topological degree (Lemma \ref{lA}), a connected subset ${\cal C}$ of zeroes joining $[0,1]\times K^-$ with $[0,1]\times K^+$ can be found (see Fig. 1). From our construction these zeroes represent geodesics joining $z$ with all the points of the  compact interval $J \subset L[x'_1,x'_2]$.

These arguments are developed rigorously in the proof of Theorem \ref{t3} and some technical properties on the $\mu_{i}$'s are postponed to a series of lemmas: (i) continuity (Lemma \ref{l5}), (ii) boundary conditions (Lemma \ref{l6}) and (iii) abstract conditions satisfied in order to apply arguments on degree (Lemma \ref{l7}).  

When conditions (\ref{e24}) are relaxed in (\ref{e28}) the new problems are, esentially: (A) perhaps the reparameterization $\gamma_i(\tau)$ naturally leads one to consider even the case $\tau = a, b$ yielding what we call a {\it fake} geodesic, and (B) conditions (\ref{lim}) may hold just for values of $K$ in a subinterval $[\bar K^-, \bar K^+] \subset [K^-,K^+]$ and, thus, the zeroes of $\mu(c,K)$ may appear in one of the four cases depicted in Fig. 2. In Section 5 we show that under our assumptions none of these possibilities are a real obstacle for geodesic connectedness; in fact, possibility (B) is related to the existence of fake geodesics, and, when these appear, equalities (\ref{e28}) allow us to provide an argument on continuity based on their escape points (Lemma \ref{l10}). Moreover, the possibility of skipping any problem when (B) happens, suggests which hypothesis may be weakened to give a more accurate result. In Section 6 we give this accurate result replacing $[\bar K^-, \bar K^+] \times [0,1]$ by regions where $K$ varies with $ c$ as depicted in Fig. 3. 

\vspace{1cm}
\begin{center}
{\bf Figs. 1, 2, 3 here}
\end{center}
\vspace{1cm}

\section{Geodesics in multiwarped spacetimes}

Let $ (F_{i},g_{i}) $ be Riemannian manifolds, $ (I,-d\tau ^{2})$  an open interval of $\R{}$ 
with $I=(a,b)$ and its usual metric reversed, and $f_{i}>0 \; i=1,\ldots ,n$  smooth functions on $I$. A {\it multiwarped spacetime} with 
base $(I,-d\tau^{2})$, fibers $(F_{i}, g_{i}) \; i=1,\ldots,n$ and warping functions $f_{i}>0$, $i=1,\ldots ,n\;$ is the product 
manifold $I\times F_{1}\times \cdots \times F_{n}$ endowed with the Lorentz metric:

\begin{equation}\label{e1}g=-\pi_{I}^{*}d\tau^{2}+\sum_{i=1}^{n}(f_{i}\circ\pi_{I})^{2}\pi^{*}_{i}g_{i}\equiv-d\tau^{2}+\sum_{i=1}^{n}f^{2}_{i}g_{i}\end{equation}
where $\pi_{I}$ and $\pi_{i} \;\; i=1,\ldots,n\;$ are the natural projections of $I\times F_{1}\times\cdots\times F_{n}$ onto I and $F_{1},\ldots,F_{n}$, respectively, and will be omitted.

 A  Riemannian manifold will be called weakly convex if any two of its points  can be joined by a geodesic which minimize the distance; if this geodesic is unique it will be called strongly convex (these names are different to those in \cite{Sa97}, see \cite{Sa98}). Denote by $d_{i}$ the distance on $F_{i}$ canonically associated to the Riemannian metric $g_{i}$. Of course, if the Riemannian manifolds $(F_{i},g_{i})$ are complete then each $F_{i}$ is weakly convex by the Hopf-Rinow theorem, but the converse is not true (see \cite{B-G-S} for a detailed study). A vector $X$ tangent to $I\times F_{1}\cdots\times F_{n}$ is lightlike if $g(X,X)=0$ and $X\neq0$, timelike if $g(X,X)<0$ and spacelike if $g(X,X)>0$ or $X=0$; the timelike vector field $\partial /\partial\tau$ fixes the canonical future orientation in $I\times F_{1}\times\cdots\times F_{n}$. Given a point $z\in I\times F_{1}\times\cdots\times F_{n}$, $I^{+}(z)$ [resp. $J^{+}(z)$] denotes the chronological [resp. causal] future of $z$ (set of points in  $I\times F_{1}\times\cdots\times F_{n}$ which can be joined by a future pointing timelike [resp. timelike or lightlike] piecewise smooth curve starting at $z$); if $z'\in I^{+}(z)$ [resp. $z'\in J^{+}(z)$] then the two points $z$, $z'$ are chronologically [resp. causally] related.

 Let $\gamma: {\cal J}\rightarrow I\times F_{1}\times\cdots\times F_{n}$, $\gamma(t)=(\tau(t),\gamma_{1}(t),\ldots,\gamma_{n}(t))$ be a (smooth) curve on the interval ${\cal J}$. Computing directly from the geodesic equations as in any warped product, it is straightforward to show that $\gamma$ is a geodesic with respect to $g$ if and only if 

\begin{equation} \label{e2} \frac{d^{2}\tau}{dt^{2}}=-\sum_{i=1}^{n}\frac{c_{i}}{f_{i}^{3}\circ\tau}\cdot\frac{df_{i}}{d\tau}\circ\tau \end{equation}

\begin{equation} \label{e3} \frac{D}{dt}\frac{d\gamma_{i}}{dt}=-\frac{2}{f_{i}\circ\tau}\cdot\frac{d(f_{i}\circ\tau)}{dt}\cdot\frac{d\gamma_{i}}{dt}\;\;\;\; i=1,\ldots ,n \end{equation}                                                                                                                                                                                                                                                                                                                                                    on ${\cal J}$, where $D/dt$ denotes the covariant derivate associated to each $g_{i}$ along $\gamma_{i}$ and $c_{i}$ is the constant $(f_{i}^{4}\circ\tau)\cdot g_{i}(d\gamma_{i}/dt,d\gamma_{i}/dt)$. Note that if $c_{i}=0$ for all $i=1,\ldots ,n$ then $d^{2}\tau /dt^{2}\equiv 0$, that is, the base of our spacetime is totally geodesic, as in any warped product. The equation (\ref{e2}) admits the first integral

\begin{equation}\label{e4}\frac{d\tau}{dt}=\epsilon\sqrt{-D+\sum_{i=1}^{n}\frac{c_{i}}{f_{i}^{2}\circ\tau}}\end{equation}
where $D=g(d\gamma /dt,d\gamma /dt)$ and $\epsilon\in \{\pm1\}$. On the other hand, by equation (\ref{e3}), each $\gamma_{i}$ is a pregeodesic of $(F_{i},g_{i})$. In fact, when $c_{i}\neq 0$ the reparametrization $\overline{\gamma}_{i}(r)=\gamma_{i}(t_{i}(r))$ where
                                                                             \begin{equation}\label{e5}\frac{dt_{i}}{dr}=\frac{1}{\sqrt{c_{i}}}f_{i}^{2}\circ\tau\circ t_{i}\end{equation}
(in a maximal domain) is a geodesic of $(F_{i},g_{i})$ being

\begin{equation}\label{e6}1=g_{i}(\frac{d\overline{\gamma}_{i}}{dr},\frac{d\overline{\gamma}_{i}}{dr}).\end{equation}.

 Let $t(\tau)$ and $r_{i}(t)$ be the inverse functions (where they exist) of the ones given by (\ref{e4}) and (\ref{e5}), respectively; then 

\begin{equation}\label{e7}\frac{d(r_{i}\circ t)}{d\tau}=\epsilon\cdot \sqrt{c_{i}} \cdot f_{i}^{-2}(-D+\sum_{j=1}^{n}\frac{c_{j}}{f_{j}^{2}})^{-1/2}\end{equation}
on a certain domain. Assume now that all the fiber components $\gamma_{i}$ of a geodesic $\gamma$ are minimizing pregeodesics and can be continuously reparametrized by $\tau$ (projection of $\gamma$ on $I$). Then, integrating (\ref{e7})    and taking into account (\ref{e2}),(\ref{e3}) and (\ref{e4}) we can find a sufficient condition for the existence of a connecting geodesic. More precisely, we obtain,

\begin{lemma}
\label{l1}
 There exists a geodesic connecting $z=(\tau_{0},x_{1},\ldots ,x_{n})$ and $z'=(\tau_{0}',x_{1}',\ldots ,x_{n}')$ if there are constants $c_{1},\ldots ,c_{n}\geq 0$, $c_{1}+\cdots +c_{n}=1$, $D\in \R{}$ such that,

(i) Either $\sum_{i=1}^{n}\frac{c_{i}}{f_{i}^{2}(\tau_{0})}\neq D$ or if this equality holds then $\frac{d}{d\tau}\sum_{i=1}^{n}\frac{c_{i}}{f_{i}^{2}}(\tau_{0})\neq 0$ and

(ii) the equality
\begin{equation}\label{e8}\int_{\tau_{0}}^{\tau_{0}'}h_{i}^{\epsilon}=l_{i}\end{equation}
 holds, where $l_{i}=d_{i}(x_{i},x_{i}')$ and the function $h_{i}^{\epsilon}(\equiv h_{i}^{\epsilon}[c_{1},\ldots ,c_{n},D]):(a_{\star},b_{\star})\subseteq I\rightarrow \R{}$, has a domain which includes $(\tau_{0},\tau_{0}')$ or $(\tau_{0}',\tau_{0})$, being equal to the right-hand side of (\ref{e7}).

\end{lemma}

 The case $\tau(t) \equiv \tau_{0}$    (which, from (\ref{e2}), is equivalent to the fact that zeroes of $\tau'(t)$ are not isolated) can be easily studied because when both equalities in formulae {\it (i)}   of Lemma 1 hold then the $\gamma_{i}'s$ are geodesics of speed $\sqrt{c_{i}}/f_{i}^{2}(\tau_{0})$, obtaining:

\begin{lemma}\label{l2}
   There exists a geodesic joining $z$ and $z'$ if $\tau_{0} =\tau_{0}'$ and there exist $c_{1},\ldots ,c_{n} \geq 0$, $c_{1}+ \cdots +c_{n}=1$ such that: $\frac{d}{d\tau}\sum_{i=1}^{n}\frac{c_{i}}{f_{i}^{2}}(\tau_{0})=0$ and                 $\frac{c_{i}}{f_{i}^{4}(\tau_{0})}=\frac{l_{i}^{2}}{f_{1}^{4}\cdot l_{1}^{2}+\cdots +f_{n}^{4}\cdot l_{n}^{2}}$ for all $i$.
\end{lemma}

Finally, the following consequence of previous formulae will be useful:

\begin{proposition} \label{p0}
If  $I\times F_{1}\times \cdots \times F_{n}$ is geodesically connected with the metric $-d\tau^2 + \sum_{i=1}^n f_i^2 g_i$ and $F_{m+1},\ldots ,F_{n}$  are strongly convex for some $m \in \{ 1,\ldots ,n-1 \}$, then $I\times F_{1}\times \cdots \times F_{m}$ is geodesically connected with $-d\tau^2 + \sum_{i=1}^m f_i^2 g_i$.
\end{proposition}

{\it Proof.} Otherwise, there are two points $(\tau_0, x_1, \dots , x_m), 
(\tau'_0, x'_1, \dots , x'_m),$
in $I\times F_1 \times \cdots \times F_m$ which cannot be joined by any geodesic; nevertheless, choosing any $x_{m+1} \in F_{m+1}, \dots , x_n \in F_n$ the points $(\tau_0, x_1, \dots , x_m , x_{m+1}, \dots x_n )$, and 
$(\tau_0', x'_1, \dots , x'_m , x_{m+1}, \dots x_n )$ can be joined by a geodesic $\gamma(t) = (\tau(t), \gamma_1(t), \dots \gamma_n(t))$ in $I\times F_1 \times \dots \times F_n$. From equality (\ref{e5}), $\gamma_{m+1}, \dots \gamma_{n}$ are reparametrizations of geodesics in $F_{m+1}, \dots , F_n$ with equal initial and final points, so, all of them are constant because of the strong convexity, which is absurd. $\diamondsuit $

\begin{remark}\label{r0}{\rm
It is possible to weaken the hypotheses on strong convexity above. In fact, it is enough for each fiber $F_i, i=m+1, \dots , n$ to admit a point $x_i$ such that no (non-constant) geodesic emanating from $x_i$ returns to $x_i$. Nevertheless, Proposition \ref{p0} does {\it not} hold if we replace strong convexity by weak convexity: inner Schwarzschild spacetime will be a nice example of this (see Theorem \ref{t6}). 
}
\end{remark}

\section{Connection by causal geodesics.} 

The following two sublemmas, even though quite obvious, are written now in order to make more readable the technical Lemma \ref{l3}.      

{\bf Sublemma 1.} {\it Consider n constants $c_{1}^{M},\ldots ,c_{n}^{M}>0$, $c_{1}^{M}+\cdots +c_{n}^{M}=k$ and fix $i \in \{1,\ldots ,n\}$. There exists arbitrarily small constants $\epsilon_{1},\ldots ,\epsilon_{i},\ldots ,\epsilon_{n}>0$ such that $c_{1}=c_{1}^{M}+\epsilon_{1},\ldots , c_{i-1}=c_{i-1}^{M}+\epsilon_{i-1}, c_{i}=c_{i}^{M}-\epsilon_{i}, 
c_{i+1}=c_{i+1}^{M}+\epsilon_{i+1},\ldots ,c_{n}=c_{n}^{M}+\epsilon_{n}$ satisfy

\begin{equation}\label{e10} \left\{ \begin{array}{l} c_{1}+ \cdots +c_{n}=k \\ \frac{c_{1}}{c_{n}}=\frac{c_{1}^{M}}{c_{n}^{M}},\ldots ,\frac{c_{i-1}}{c_{n}}=\frac{c_{i-1}^{M}}{c_{n}^{M}},\frac{c_{i+1}}{c_{n}}=\frac{c_{i+1}^{M}}{c_{n}^{M}},\ldots ,\frac{c_{n-1}}{c_{n}}=\frac{c_{n-1}^{M}}{c_{n}^{M}} \end{array} \right. \end{equation}} 

{\it Proof.} Choose $\epsilon_{n}=\epsilon <k-c_{n}^{M}$, $\epsilon_{j}=\frac{c_{j}^{M}}{c_{n}^{M}}\epsilon$ for $j\neq i,n$ and $\epsilon_{i}=\sum_{j\neq i}\epsilon_{j}$. $\diamondsuit $

{\bf Sublemma 2.} {\it If $c_{1},\ldots ,c_{n},\overline{c}_{1},\ldots ,\overline{c}_{n}>0$, $c_{1}+\cdots +c_{n} = 1 = \overline{c}_{1}+\cdots +\overline{c}_{n}$ and $\overline{c}_{n}<c_{n}$, then there exists $i_{0}\in \{1,\ldots ,n-1\}$ such that

\begin{equation}\label{e11}\begin{array}{l}\ \overline{c}_{i_{0}}>c_{i_{0}} \quad    (and, thus, \frac{\overline{c}_{n}}{\overline{c}_{i_{0}}}<\frac{c_{n}}{c_{i_{0}}}) \\ 
   \frac{\overline{c}_{j}}{\overline{c}_{i_{0}}} \leq \frac{c_{j}}{c_{i_{0}}} \quad \forall j\neq n \end{array}\end{equation} 

Moreover, if there also exists $j\in \{1,\ldots ,n-1\}$ such that $\frac{\overline{c}_{j}}{\overline{c}_{n}}<\frac{c_{j}}{c_{n}}$ then there exists $j_{0}\in \{1,\ldots ,n-1\}$ such that 

\begin{equation}\label{e12}\begin{array}{l}\ \overline{c}_{j_{0}}<c_{j_{0}},\quad  \frac{\overline{c}_{n}}{\overline{c}_{j_{0}}}>\frac{c_{n}}{c_{j_{0}}} \\  \frac{\overline{c}_{j}}{\overline{c}_{j_{0}}} \geq \frac{c_{j}}{c_{j_{0}}}\end{array}\end{equation}}

{\it Proof.} For the first asertion, take $i_{0}\in \{1,\ldots ,n-1\}$ such that $\frac{\overline{c}_{i_{0}}}{c_{i_{0}}}$ is a maximum, and for the last one take $\frac{\overline{c}_{j_{0}}}{c_{j_{0}}}$ minimum. $\diamondsuit $

The following result will be essential to reduce our problem.

\begin{lemma}
\label{l3}
Let $(\tau_{0},\tau'_{0})$ be an interval of $\;\R{}$, and $ w_{1}, \ldots ,w_{n}\;$ n (smooth) positive functions defined on $(\tau_{0},\tau'_{0})$. Take constants $c_{1}, \ldots ,c_{n}\geq 0$, $c_{1}+ \cdots +c_{n}=k$, $0<k\leq 1$ and $l_{1}, \ldots ,l_{n}\geq 0$, $ l_{1}+ \cdots +l_{n}>0$. Let 
$\nu\geq 0$ be a smooth function. If

\begin{equation} \label{e13} \left\{ \begin{array}{c} \sqrt{c_{1}}\int_{\tau_{0}}^{\tau'_{0}}w_{1}(c_{1}w_{1}+ \cdots +c_{n}w_{n}+\nu)^{-1/2} \geq l_{1} \\ \vdots \\ \sqrt{c_{n}}\int_{\tau_{0}}^{\tau'_{0}}w_{n}(c_{1}w_{1}+ \cdots +c_{n}w_{n}+\nu)^{-1/2} \geq l_{n} \end{array} \right. \end{equation} 
then there exists unique   $c'_{1},\cdots,c'_{n}\geq 0$, $c'_{1}+\cdots +c'_{n}=k$ and $D\leq 0$, which vary continuously with $k$ and $\nu $ (in the $L^{\infty}[\tau_{0},\tau'_{0}]$ topology), such that

\begin {equation}\label{e14} \left\{ \begin{array}{c} \sqrt{c'_{1}}\int_{\tau_{0}}^{\tau'_{0}}w_{1}(c'_{1}w_{1}+ \cdots +c'_{n}w_{n}+\nu-D)^{-1/2}= l_{1} \\ \vdots \\ \sqrt{c'_{n}}\int_{\tau_{0}}^{\tau'_{0}}w_{n}(c'_{1}w_{1}+ \cdots +c'_{n}w_{n}+\nu-D)^{-1/2}= l_{n} \end{array} \right. \end{equation}

\end{lemma}

{\it Proof}. We will work by induction on the number $n$ of functions. The case $n=1$ is obvious, and assume it is true for $n-1$. Consider the compact set $S=\{(c_{1},\ldots,c_{n})\in \R{}^{n}: c_{1},\ldots,c_{n}\geq 0, c_{1}+\cdots +c_{n}=k,\quad \hbox{satisfying} \; (\ref{e13}) \}$. By hypothesis, $S\neq \emptyset$ and we can assume all the $c_{i}'s$ are non null (otherwise, some $l_{i}$ is $0$ and the result is trivial from the induction hypothesis). Let $(c_{1}^{M},\ldots ,c_{n}^{M})\in S$ be such that
 
\begin{equation}\label{e15} \sqrt{c_{n}^{M}}\int_{\tau_{0}}^{\tau'_{0}}w_{n}(c_{1}^{M}w_{1}+\cdots +c_{n}^{M}w_{n}+\nu)^{-1/2}\end{equation}
is a maximum on S. Then the equality in (\ref{e13}) must hold for the $n-1$ first inequalities. In fact, if there exists $i \in \{1,\ldots ,n-1\}$ which does not satisfy the equality, then we can choose $c_{1},\ldots ,c_{n}$ as in Sublemma 1. For $\epsilon_{1},\ldots ,\epsilon_{n}$ small enough, the i-th inequality (\ref{e13}) holds and dividing all the equations (\ref{e10}) by $\frac{c_{j}}{c_{n}}=\frac{c_{j}^{M}}{c_{n}^{M}}$ it is clear that (\ref{e13}) remains true for $j=1,\ldots ,i-1,i+1,\ldots ,n-1$. Now, equations (\ref{e10}) not only shows that the last equation (\ref{e13}) holds, but also clearly contradicts the maximality assumption on (\ref{e15}).  

For each $c_{n}\in [0,c_{n}^{M}]\;$ we can put $k'=k-c_{n}$, $\;\nu'=\nu+c_{n}w_{n}$ and obtain unique $(c'_{1},\ldots ,c'_{n-1})\equiv (c_{1}(c_{n}),\ldots ,c_{n-1}(c_{n}))$, $\;D\equiv D(c_{n})\;$ depending continuously on $c_{n}$ by using the induction hypothesis. If $c_{n}=c_{n}^{M}$ then $D(c_{n}^{M})=0$, $c_{1}(c_{n}^{M})=c_{1}^{M},\ldots ,c_{n-1}(c_{n}^{M})=c_{n-1}^{M}$, thus if we define $G_{n}:[0,c_{n}^{M}]\rightarrow \R{}$

\[G_{n}(c_{n})=\sqrt{c_{n}}\int_{\tau_{0}}^{\tau'_{0}}w_{n}(c_{1}(c_{n})w_{1}+\cdots +c_{n-1}(c_{n})w_{n-1}+c_{n}w_{n}+\nu -D(c_{n}))^{-1/2}-l_{n}\]
we have $G_{n}(c_{n}^{M})\geq 0$. Clearly, $G_{n}(0)\leq 0$ and if $G_{n}$ vanishes at $c'_{n}\in [0,c_{n}^{M}]$ then the corresponding $(c_{1}(c'_{n}),\ldots ,c_{n-1}(c'_{n}),c'_{n})$ and $D(c'_{n})$ satisfy (\ref{e14}). As $G_{n}$ is continuous, by the induction hypothesis, and it varies continuously with $k$ and $\nu $, it is necessary to prove just the uniquess of $c'_{n}$, which will be checked by showing that $G_{n}$ is strictly increasing. Assume $\overline{c}_{n}<c_{n}$ and put $c_{1}\equiv c_{1}(c_{n}),\ldots ,\overline{c}_{n-1}\equiv c_{n-1}(\overline{c}_{n})$. Let $i_{0}\in \{1,\ldots ,n-1 \}$ be given by Sublemma 2. As $D(\overline{c}_{n})$, $D(c_{n})$ satisfy

\begin{equation}\label{e16} \begin{array}{l} \quad \quad l_{i_{0}}= \\ \int_{\tau_{0}}^{\tau'_{0}}w_{i_{0}}\left(\frac{c_{1}}{c_{i_{0}}}w_{1}+\cdots +\frac{c_{i_{0}-1}}{c_{i_{0}}}w_{i_{0}-1}+w_{i_{0}}+\frac{c_{i_{0}+1}}{c_{i_{0}}}w_{i_{0}+1}+\cdots +\frac{c_{n}}{c_{i_{0}}}w_{n}+\frac{\nu}{c_{i_{0}}}-\frac{D(c_{n})}{c_{i_{0}}}\right)^{-1/2}= \\ \int_{\tau_{0}}^{\tau'_{0}}w_{i_{0}}\left(\frac{\overline{c}_{1}}{\overline{c}_{i_{0}}}w_{1}+\cdots +\frac{\overline{c}_{i_{0}-1}}                                                                                                                                             
{\overline{c}_{i_{0}}}w_{i_{0}-1}                                                                                                                                                                                              +w_{i_{o}}+\frac{\overline{c}_{i_{0}+1}}{\overline{c}_{i_{0}}}w_{i_{0}+1}+\cdots +\frac{\overline{c}_{n}}{\overline{c}_{i_{0}}}w_{n}+\frac{\nu}{\overline{c}_{i_{0}}}-\frac{D(\overline{c}_{n})}{\overline{c}_{i_{0}}}\right)^{-1/2} \end{array} \end{equation} 
thus, $ - \frac{D(\overline{c}_{n})}{\overline{c}_{i_{0}}}>- \frac{D(c_{n})}{c_{i_{0}}} $. This inequality and (\ref{e11}) imply

\begin{equation}\label{e17} - \frac{D(\overline{c}_{n})}{\overline{c}_{n}}>- \frac{D(c_{n})}{c_{n}}, \end{equation} thus if 

\begin{equation}\label{e18}\frac{\overline{c}_{j}}{\overline{c}_{n}} \geq \frac{c_{j}}{c_{n}}\end{equation} for all $j$, it is clear that $G_{n}(\overline{c}_{n})<G_{n}(c_{n})$. If (\ref{e18}) does not hold then taking $j_{0}\in \{1,\ldots ,n-1\}$ as in the last assertion of Sublemma 2, and using an inequality analogous to (\ref{e17}) for $j_{0}$: 

\begin{equation}\label{e19} \begin{array}{c} \frac{c_{1}}{c_{j_{0}}}w_{1}+\cdots +\frac{c_{j_{0}-1}}{c_{j_{0}}}w_{j_{0}-1}+w_{j_{0}}+\frac{c_{j_{0}+1}}{c_{j_{0}}}w_{j_{0}+1}+\cdots +\frac{c_{n}}{c_{j_{0}}}w_{n}+\frac{\nu}{c_{j_{0}}}-\frac{D(c_{n})}{c_{j_{0}}}< \\ \frac{\overline{c}_{1}}{\overline{c}_{j_{0}}}w_{1}+\cdots +\frac{\overline{c}_{j_{0}-1}}                                                                                                                                             
{\overline{c}_{j_{0}}}w_{j_{0}-1}                                                                                                                                                                                              +w_{j_{o}}+\frac{\overline{c}_{j_{0}+1}}{\overline{c}_{j_{0}}}w_{j_{0}+1}+\cdots +\frac{\overline{c}_{n}}{\overline{c}_{j_{0}}}w_{n}+\frac{\nu}{\overline{c}_{j_{0}}}-\frac{D(\overline{c}_{n})}{\overline{c}_{j_{0}}}. \end{array} \end{equation}

But this is a contradiction with

\begin{equation}\label{e20} \begin{array}{l} \quad \quad l_{j_{0}}= \\ \int_{\tau_{0}}^{\tau'_{0}}w_{j_{0}}\left(\frac{c_{1}}{c_{j_{0}}}w_{1}+\cdots +\frac{c_{j_{0}-1}}{c_{j_{0}}}w_{j_{0}-1}+w_{j_{0}}+\frac{c_{j_{0}+1}}{c_{j_{0}}}w_{j_{0}+1}+\cdots +\frac{c_{n}}{c_{j_{0}}}w_{n}+\frac{\nu}{c_{j_{0}}}-\frac{D(c_{n})}{c_{j_{0}}}\right)^{-1/2}= \\ \int_{\tau_{0}}^{\tau'_{0}}w_{j_{0}}\left(\frac{\overline{c}_{1}}{\overline{c}_{j_{0}}}w_{1}+\cdots +\frac{\overline{c}_{j_{0}-1}}                                                                                                                                             
{\overline{c}_{j_{0}}}w_{j_{0}-1}                                                                                                                                                                                              +w_{j_{o}}+\frac{\overline{c}_{j_{0}+1}}{\overline{c}_{j_{0}}}w_{j_{0}+1}+\cdots +\frac{\overline{c}_{n}}{\overline{c}_{j_{0}}}w_{n}+\frac{\nu}{\overline{c}_{j_{0}}}-\frac{D(\overline{c}_{n})}{\overline{c}_{j_{0}}}\right)^{-1/2}. \end{array} \end{equation} $\diamondsuit $

                                                                                                                                                                           Lemma \ref{l3} allows us to obtain the following result on multiwarped spacetimes.
\begin{theorem}
\label{t2}
Let $(I\times F_{1} \times \cdots \times F_{n}  , g)$ be a multiwarped   spacetime  with  weakly convex 
fibers, and consider any pair of distinct points 
$(\tau_{0}, x_{1},\ldots , x_{n}), (\tau_{0}',x_{1}',\ldots , x_{n} ') \in  I \times F_{1} \times \cdots \times F_{n}  , \tau_{0} \leq \tau_{0}'$. 

     (1) The following conditions are equivalent:
\begin{quote}
          (i) There exists a timelike geodesic joining  $(\tau_{0}, x_{1},\ldots , x_{n})$ and $ (\tau_{0}',x_{1}',\ldots , x_{n} ').$

(ii) There exists $c_{1},\ldots ,c_{n} \geq 0$, $c_{1}+\cdots +c_{n}=1$ such that 

\begin{equation} \label{e21} \sqrt{c_{i}}\int_{\tau_{0}}^{\tau'_{0}}f_{i}^{-2}(\frac{c_{1}}{f_{1}^{2}}+ \cdots +\frac{c_{n}}{f_{n}^{2}})^{-1/2} \geq l_{i} \end{equation} 
for all $i$, where $l_{i}=d_{i}(x_{i},x_{i}')$ and with equality in j-th equation if and only if $c_{j}=0$. 
        
(iii) $(\tau_{0}',x_{1}',\ldots , x_{n} ') \in I^{+}(\tau_{0}, x_{1},\ldots , x_{n}).$
\end{quote}

     (2) The following conditions are also equivalent:
\begin{quote}
(i) There exists a timelike or lightlike geodesic joining 
$(\tau_{0}, x_{1},\ldots , x_{n})$ and $ (\tau_{0}',x_{1}',\ldots , x_{n} ') $.

(ii) There exists $c_{1},\ldots ,c_{n} \geq 0$, $c_{1}+\cdots +c_{n}=1$ such that inequalities (\ref{e21}) hold.

(iii) $(\tau_{0}',x_{1}',\ldots , x_{n} ') \in J^{+}(\tau_{0}, x_{1},\ldots ,x_{n}).$
\end{quote}

      Moreover, if the equality holds in all equations (\ref{e21}) then there is a lightlike and no timelike geodesic joining the points.

\end{theorem}

{\it Proof.} Implications (i)$\Rightarrow$ (iii) are obvious; for the remainder use Lemma \ref{l3} and a reasoning as in \cite[Theorems 3.7, 4.2]{Sa97}. $\diamondsuit$

\section{Geodesic connectedness. Special case.}

As a first claim, Lemma \ref{l1} can be extended to cases where $\tau'(t)$ has isolated zeroes (and, thus, $\tau$ can be used as a parameter around all points but these zeroes). In fact, note that the zeroes of $\tau'$ are isolated if and only if last inequality in Lemma \ref{l1}(i) holds at each zero and, in this case, the zero is a strict relative maximum or minimum. More precisely, fix $c_{1},\ldots ,c_{n}\geq 0$, $c_{1}+\cdots +c_{n}=1$, $D\in \R{}$ such that $\frac{c_{1}}{f_{1}^{2}(\tau_{0})}+\cdots +\frac{c_{n}}{f_{n}^{2}(\tau_{0})}\geq D$ and consider the subsets
\begin{equation}\label{e22}A_{+}=\{\tau\in (a,b): \tau_{0}\leq \tau, \sum_{i=1}^{n}\frac{c_{i}}{f_{i}^{2}(\tau)}=D\}\cup\{b\}\end{equation}
\begin{equation}\label{e23}A_{-}=\{\tau\in (a,b): \tau_{0}\geq \tau, \sum_{i=1}^{n}\frac{c_{i}}{f_{i}^{2}(\tau)}=D\}\cup\{a\}\end{equation}
 Define $a_{\star}\equiv a_{\star}(c_{1},\ldots ,c_{n},D)$, $b_{\star}\equiv b_{\star}(c_{1},\ldots ,c_{n},D)$ by:

\begin{equation}\label{e23.5} \hbox{If} \; \frac{d}{d\tau}\sum_{i=1}^{n}\frac{c_{i}}{f_{i}^{2}}(\tau_{0})\left \{ \begin{array}{c} >0 \\ <0 \\ =0 \end{array} \right. \; \hbox{then} \; b_{\star}= \left \{ \begin{array}{c}  min(A_{+}-\{\tau_{0}\})  \\ min(A_{+}), \\ min(A_{+}) \end{array} \right.  a_{\star}= \left \{ \begin{array}{c} max(A_{-}) \\ max(A_{-}-\{\tau_{0}\}) \\ max(A_{-}) \end{array} \right. \end{equation}
  
 Lemma \ref{l1} also holds if we assume the following convention for (\ref{e8}).

 \begin{convention}\label{c1} {\rm From now on integral (\ref{e8}) will be understood in the following generalized sense: if $\int_{\tau_{0}}^{b_{\star}}h_{i}^{\epsilon=1}\geq l_{i}$ for some $i$, then the first member of (\ref{e8}) denotes the usual integral; otherwise and if $b_{\star}\neq b$, we can follow integrating, by reversing the sense of integration and, if $\int_{\tau_{0}}^{b_{\star}}h_{i}^{\epsilon=1}-\int_{b_{\star}}^{a_{\star}}h_{i}^{\epsilon=1}\geq l_{i}$ for some $i$ then the first member of (\ref{e8}) means $\int_{\tau_{0}}^{b_{\star}}h_{i}^{\epsilon=1}-\int_{b_{\star}}^{\tau_{0}'}h_{i}^{\epsilon=1}$ for all $i$. If this last inequality does not hold for any $i$ and $a_{\star}\neq a$, the procedure must follow reversing the sense of integration as many times as necessary in the obvious way. Analogously, when $\epsilon=-1$, condition (\ref{e8}) means either $\int_{\tau_{0}}^{\tau_{0}'}h_{i}^{\epsilon=-1}$ (for all $i$; in this case necessarily $\tau_{0}'< \tau_{0}$) or    $\int_{\tau_{0}}^{a_{\star}}h_{i}^{\epsilon=-1}-\int_{a_{\star}}^{\tau_{0}'}h_{i}^{\epsilon=-1}$ (for all $i$) or $\int_{\tau_{0}}^{a_{\star}}h_{i}^{\epsilon=-1}-\int_{a_{\star}}^{b_{\star}}h_{i}^{\epsilon=-1}+\int_{b_{\star}}^{\tau_{0}'}h_{i}^{\epsilon=-1}$ (for all $i$), etc.}
\end{convention}

\begin{theorem}\label{t3} A spacetime $(I\times F_{1}\times \cdots \times F_{n},g)$ with weakly convex fibers $(F_{1},g_{1}),\ldots ,$ $(F_{n},g_{n})$ satisfying
\begin{equation} \label{e24}
\begin{array}{ll}
\int_{c}^{b}f_{i}^{-2}(\frac{1}{f_{1}^{2}}+\cdots +\frac{1}{f_{n}^{2}}+1)^{-1/2}=\infty  &  \int_{a}^{c}f_{i}^{-2}(\frac{1}{f_{1}^{2}}+\cdots +\frac{1}{f_{n}^{2}}+1)^{-1/2}=\infty \end{array} 
\end{equation}
for all $i\in \{1,\ldots ,n\}$ and for $c\in (a,b)$, is geodesically connected.
\end{theorem}

The hypothesis of Theorem \ref{t3} will be useful in the sense of the following Lemma.

\begin{lemma}
\label{l4}
Equalities (\ref{e24}) imply 
\begin{equation}\label{e25} \begin{array}{ll} \int_{c}^{b}f_{i}^{-2}(\frac{c_{1}}{f_{1}^{2}}+\cdots +\frac{c_{n}}{f_{n}^{2}}-D)^{-1/2}=\infty  & \int_{a}^{c'}f_{i}^{-2}(\frac{c_{1}}{f_{1}^{2}}+\cdots +\frac{c_{n}}{f_{n}^{2}}-D)^{-1/2}=\infty \end{array} \end{equation}
for all $i$ and for all $c_{1},\ldots ,c_{n}\geq 0$, $c_{1}+\cdots +c_{n}=1$ and $c,c'\in I$, $D\in \R{}$ such that the denominators of the integrals do not vanish.
\end{lemma}

{\it Proof of Lemma \ref{l4}.} Since $\int_{c}^{b}f_{i}^{-2}(\frac{1}{f_{1}^{2}}+\cdots +\frac{1}{f_{n}^{2}}+1)^{-1/2}=\infty$, then  multiplying the integrand by the function $\rho(\tau)$ such that   $\rho(\tau)\cdot (\frac{1}{f_{1}^{2}}+\cdots +\frac{1}{f_{n}^{2}}-D)^{1/2}=(\frac{1}{f_{1}^{2}}+\cdots +\frac{1}{f_{n}^{2}}+1)^{1/2}$ (which satisfies $\rho(\tau)\geq \delta_{0}>0$) we would have that                                        $\int_{c}^{b}f_{i}^{-2}(\frac{1}{f_{1}^{2}}+\cdots +\frac{1}{f_{n}^{2}}-D)^{-1/2}=\infty$. As  $(\frac{c_{1}}{f_{1}^{2}}+\cdots +\frac{c_{n}}{f_{n}^{2}}-D)^{-1/2}\geq c_{k}^{-1/2}\cdot (\frac{1}{f_{1}^{2}}+\cdots +\frac{1}{f_{n}^{2}}-\frac{D}{c_{k}})^{-1/2}$ where $c_{k}$ is the maximum of the $c_{i}'s$, the result follows obviously. $\diamondsuit$

In what follows $\hat{c}=(c_{1},\ldots ,c_{n})$.

{\it Proof of Theorem \ref{t3}.} Let $(\tau_{0},x_{1},\ldots ,x_{n}), \; (\tau_{0}',x_{1}',\ldots ,x_{n}') \in I\times F_{1}\times\cdots \times F_{n}$  and put $l_{i}=d_{i}(x_{i},x_{i}')$. We can assume $l_{i}>0$ for all $i$ because, otherwise, the problem would be reduced to the case of $n-1$ fibers just by putting $\gamma_{i}\equiv constant$ for the corresponding $i$.

Put $K^{+}=max_{i}\{\frac{1}{f_{i}^{2}(\tau_{0})}\}$ and $K^{-}=-K^{+}$, and  define the set $\bigtriangleup_{n}=\{\hat{c}\in (0,1)^{n}: c_{1}+\cdots +c_{n}=1\}$. Fixed $(\hat{c},D)\in \bigtriangleup_{n}\times \R{}$ with $D\leq \frac{c_{1}}{f_{1}^{2}(\tau_{0})}+\cdots +\frac{c_{n}}{f_{n}^{2}(\tau_{0})}$, we consider the functions $h_{i}^{\epsilon}$ given in Lemma \ref{l1} and $a_{\star}, b_{\star}$ in (\ref{e23.5}). Then define the functions $s_{i}^{+}:\bigtriangleup_{n}\times ]0,K^{+}]\rightarrow \R{}$, $s_{i}^{+}(\hat{c},K)=t_{i}-\tau_{0}$ where $\int_{\tau_{0}}^{t_{i}}h_{i}^{\epsilon=1}[\hat{c},D]=l_{i}$ and $D=\frac{c_{1}}{f_{1}^{2}(\tau_{0})}+\cdots +\frac{c_{n}}{f_{n}^{2}(\tau_{0})}-K$, taking into account Convention \ref{c1} for each $i$ that is, if $\int_{\tau_{0}}^{b_{\star}}h_{i}^{\epsilon=1}-\int_{b_{\star}}^{t_{i}}h_{i}^{\epsilon=1}=l_{i}$ then, actually,  $s_{i}^{+}(\hat{c},K)=2b_{\star}-t_{i}-\tau_{0}$, etc. (note that we can use Convention \ref{c1}, since $b_{\star}\neq b$ and $a_{\star}\neq a$ because of Lemma \ref{l4}). We define $s_{i}^{-}:\bigtriangleup_{n}\times [K^{-},0[\rightarrow \R{}$ analogously but using $h_{i}^{\epsilon=-1}[\hat{c},D]$ and opposite sense of integration, and being $D=\frac{c_{1}}{f_{1}^{2}(\tau_{0})}+\cdots +\frac{c_{n}}{f_{n}^{2}(\tau_{0})}+K$, that is, $s_{i}^{-}(\hat{c},K)=\tau_{0}-t_{i}$ where $\int_{\tau_{0}}^{t_{i}}h_{i}^{\epsilon=-1}=l_{i}$ (if $\int_{\tau_{0}}^{a_{\star}}h_{i}^{\epsilon=-1}+\int_{t_{i}}^{a_{\star}}h_{i}^{\epsilon=-1}=l_{i}$ then $s_{i}^{-}=\tau_{0}+t_{i}-2a_{\star}$, etc.). From Lemma \ref{l5} each couple of functions $s_{i}^{+}$, $s_{i}^{-}$ glues together in a continuous function $s_{i}$ on all $\bigtriangleup_{n}\times [K^{-},K^{+}]$; moreover, also from Lemma \ref{l5} the functions $\mu_{i}=1-\frac{s_{i}}{s_{1}}\quad i=1,\ldots ,n$ are continuously well defined on all $\bigtriangleup_{n}\times [K^{-},K^{+}]$. Recall that for the geodesic determined by $\hat{c},D$  the value of $K$ corresponds to $sign(\tau'(t_{0}))\cdot \tau'(t_{0})^{2}$ thus, our choice of $K^{+},K^{-}$ implies that any geodesic characterized by $(\hat{c},K^{\pm})$ is causal.                
 
Consider the homeomorphism ${\cal Y}:\bigtriangleup_{n}\rightarrow (0,1)^{n-1}$,
\begin{equation}\label{e25.5} {\cal Y}(\hat{c})=(\frac{c_{2}}{c_{1}+c_{2}},\ldots ,\frac{c_{n-1}}{c_{1}+\cdots +c_{n-1}},c_{n})\equiv (y_{1},\ldots ,y_{n-1}) \equiv \hat{y}\end{equation}
 and define $\bar{\mu_{i}}:(0,1)^{n-1}\times [K^{-},K^{+}]\rightarrow \R{}$, $\bar{\mu_{i}}(\hat{y},K)=\mu_{i}({\cal Y}^{-1}(\hat{y}),K)$ (in order to simplify the notation, we also denote $\bar{\mu_{i}}(\hat{y},K)$ by $\mu_{i}(\hat{y},K)$ and, in general, ${\cal Y}$ will be omitted if there is no possibility of confusion). Next, we will prove the following essential step.

\begin{lemma} \label{lA}  The set of points $(\hat{y},K) \in (0,1)^{n-1}\times [K^{-},K^{+}]$ verifying $\mu_{2}(\hat{y},K)=\cdots =\mu_{n}(\hat{y},K)=0$ admits a connected component $\cal{C}$ such that

\begin{equation}\label{e25.75} {\cal C}\cap \{(0,1)^{n-1}\times \{K^{+}\}\}\neq\emptyset \; \hbox{and} \; {\cal C}\cap \{(0,1)^{n-1}\times \{K^{-}\}\}\neq \emptyset \end{equation} 
\end{lemma}

{\it Proof of Lemma \ref{lA}.} From Lemma \ref{l6}, fixed $\delta >0$, there exists $\epsilon >0$ such that $\mu_{i}(\hat{c},K)<\mu_{j}(\hat{c},K)$ being $c_{i}\leq \epsilon$, $c_{j}\geq \delta$ and $K\in [K^{-},K^{+}]$. Now apply Lemma \ref{l7}, to obtain the corresponding functions $\tilde{\mu_{i}} \; \; i=1,\ldots ,n$. From Lemma \ref{l7}(i), we have just to prove that the set of points $(\hat{y},K)\in [-1,2]^{n-1}\times [K^{-},K^{+}]$ verifying $\tilde{\mu_{2}}(\hat{y},K)=\cdots =\tilde{\mu_{n}}(\hat{y},K)=0$ admits a connected component $\hat{\cal{C}}$ such that $\hat{\cal {C}}\cap \{[-1,2]^{n-1}\times \{K^{+}\}\}\neq\emptyset$ and $\hat{\cal {C}}\cap \{[-1,2]^{n-1}\times \{K^{-}\}\}\neq \emptyset$. 
 
Now, define the function
\[
\begin{array}{llll}
 {\cal F}: & [-1,2]^{n-1}\times [K^{-},K^{+}] & \rightarrow & \R{}^{n-1} \\

           & (y_{1},\ldots ,y_{n-1},K)       & \rightarrow & (\tilde{\mu_{2}},\ldots ,\tilde{\mu_{n}})(y_{1},\ldots ,y_{n-1},K) \\    &  &   & +(y_{1},\ldots ,y_{n-1})
\end{array}
\]
and, for each $K$, let ${\cal F}_K$ be  the function obtained by fixing the last variable of ${\cal F}$ equal to $K$. Let  $\hat{\cal{F}}$ be the affine map constructed in {\it (iii)} of 
Lemma \ref{l7} plus the identity; clearly, at each point of its range  we have deg(${\cal F}_K$)= deg($\hat{{\cal F }}$) =1  (the map $\lambda \mapsto {\cal F}_K+\lambda (\hat{{\cal F}}-{\cal F}_K)$, $\lambda \in [0,1]$ is a homotopy from ${\cal F}_K$ to $\hat{{\cal F}}$ which fixes the boundary).  
Thus, if the number of solutions to the equation in $\hat y$
\begin{equation}
\label{LS}
\hat y - {\cal F}(\hat y,K) = 0
\end{equation}
were finite for some value of $K$ then the sum of their degrees would be equal to 1 and, according to \cite[Theor. Fond.]{L-S}, the required connected component of solutions to (\ref{LS}) would be found. This hypothesis of finiteness, in principle, may not hold, even though it can be removed (see \cite[Proposition 2.3]{FLN} and \cite[Lemma 3.4]{N}); anyway, the following simple reasoning skip this problem in our case, and will be convenient later. For each $m\in \N{}$ big enough, put $K^-_m = K^- + 1/m$ and construct a perturbation ${\cal F}^m$ of ${\cal F}$ satisfying: (i) ${\cal F}^m_K = {\cal F}_K$ for $K\in [K^-_m,K^+]$, (ii) ${\cal F}^m_K = \hat{{\cal F}}$ for $K =  K^-$, (iii) ${\cal F}^m_K$ is an homotopy between $\hat{\cal F}$ and ${\cal F}_{K^-_m}$  for $K\in [K^-,K^-_m]$. Now, there are just a solution to equation (\ref{LS}) for ${\cal F}^m $ when $K=K_-$, and, then, a connected component $\hat{{\cal C}}^m$ of solutions satisfying (\ref{e25.75}) is found. Taking ${\cal C}^m$ as a connected component of 
$\hat{{\cal C}}^m-([-1,2]^{n-1}\times [K^-,K^-_m))$ satisfying the first inequation in (\ref{e25.75}) then ${\cal C}^m$ is a set of solutions of the original equation (\ref{LS}) for the non-perturbed ${\cal F}$, and, by standard arguments (see \cite[Ch.I Theor. 9.1]{Wh}), limsup$\{ {\cal C}^m \}$ is the required connected component.
$\diamondsuit $

From the construction of $\cal{C}$, Lemma \ref{l2} and Lemma \ref{l5}, every point of this set represents a geodesic joining $z$ and the line $L[x_{1}',\ldots ,x_{n}']$. 

As the timelike component of a geodesic depends on $\hat{c}$ and $\tau'(t_{0})$ continuously, thus it depends on $\hat{y}$ and $K$ continuously too. Since $\cal{C}$  is connected, we obtain that $z=(\tau_{0},x_{1},\ldots ,x_{n})$ can be connected with an interval $J$ of $L[x_{1}',\ldots ,x_{n}']$ by means of geodesics. Moreover, as $\cal{C}$ satisfies (\ref{e25.75}), the extremes of the interval $J$ are one in $J^{+}(z)$ and one in $J^{-}(z)$. If $z'$ belongs to interval $J$ the proof is over, otherwise, $z'\in J^{+}(z)\cup J^{-}(z)$ and Theorem \ref{t2} can be claimed. $\diamondsuit$

\begin{lemma}
\label{l5}
For each $i=1,\ldots ,n$, there is a (unique) continuous function $s_{i}:\bigtriangleup_{n}\times [K^{-},K^{+}]\rightarrow \R{}$ whose restriction to $\bigtriangleup_{n}\times ]0,K^{+}]$ (resp.  $\bigtriangleup_{n}\times [K^{-},0[$) agrees $s_{i}^{+}$ (resp. $s_{i}^{-}$). Moreover, the functions  
$\mu_{i}=1-\frac{s_{i}}{s_{1}}\quad i=1,\ldots ,n$ ($\mu_{1}\equiv 0$) defined where makes sense admit continuous extensions to $\bigtriangleup_{n}\times [K^{-},K^{+}]$ which satisfy $\mu_{i}(\hat{c},K)=1-\frac{\sqrt{c_{1}}\cdot l_{i}\cdot f_{i}^{2}(\tau_{0})}{\sqrt{c_{i}}\cdot l_{1}\cdot f_{1}^{2}(\tau_{0})}$             when $s_{1}(\hat{c},K)=0$.                
\end{lemma}

{\it Proof.} Previously, consider each $s_{i}^{\pm}$ as a function of $(\hat{c},D)$. In order to prove the assertion on the $s_{i}$'s , we will check that every convergent sequence $\{(\hat{c}^{k},D^{k})\}_{k\in \N{}}$, $(\hat{c}^{k},D^{k})\rightarrow (\hat{c}^{\infty},D^{\infty})$, with $\hat{c}^{\infty}\in \bigtriangleup_{n}$, satisfies  $s_{i}^{+}(\hat{c}^{k},D^{k})\rightarrow s_{i}^{+}(\hat{c}^{\infty},D^{\infty})$ and $s_{i}^{-}(\hat{c}^{k},D^{k})\rightarrow s_{i}^{-}(\hat{c}^{\infty},D^{\infty})$ where makes sense and, when $s_{i}^{+}$ and $s_{i}^{-}$ are not defined on $(\hat{c}^{\infty},D^{\infty})$, $lim_{k\rightarrow \infty}s_{i}^{+}(\hat{c}^{k},D^{k})=lim_{k\rightarrow \infty}s_{i}^{-}(\hat{c}^{k},D^{k})$.

Consider the corresponding  $a_{\star}^{k}\equiv a_{\star}^{k}(\hat{c}^{k},D^{k})$, $b_{\star}^{k}\equiv b_{\star}^{k}(\hat{c}^{k},D^{k})$,    $a_{\star}^{\infty}$,$b_{\star}^{\infty}$ defined in (\ref{e23.5}); recall that  $[a_{\star}^{\infty},b_{\star}^{\infty}]\subseteq [lim sup \; a_{\star}^{k}, lim inf \; b_{\star}^{k}]$. We can consider the following cases: 

(i) If $\frac{d}{d\tau}\sum_{i=1}^{n}\frac{c_{i}^{\infty}}{f_{i}^{2}}\mid_{a_{\star}^{\infty},b_{\star}^{\infty}}\neq 0$ then the sequence of intervals $(a_{\star}^{k},b_{\star}^{k})$ converges to $(a_{\star}^{\infty},b_{\star}^{\infty})$ and $\int_{a_{\star}^{\infty}}^{b_{\star}^{\infty}}h_{i}^{\epsilon=1}[\hat{c}^{\infty},D^{\infty}]<\infty$. Moreover, $\{\sum_{i=1}^{n}\frac{c_{i}^{k}}{f_{i}^{2}}\}_{k\in \N{}}$   and its first derivatives converge uniformly  to $\sum_{i=1}^{n}\frac{c_{i}^{\infty}}{f_{i}^{2}}$   and its first derivative on an open interval containing $[a_{\star}^{\infty},b_{\star}^{\infty}]$, then the functions $h_{i}^{\epsilon}[\hat{c}^{k},D^{k}]$ converge almost everywhere to $h_{i}^{\epsilon}[\hat{c}^{\infty},D^{\infty}]$ and are dominated by an integrable function on $[a_{\star}^{\infty},b_{\star}^{\infty}]$  which implies the required convergence.

(ii) If either $b_{\star}^{\infty}=b$, $a_{\star}^{\infty}=a$, or $\frac{d}{d\tau}\sum_{i=1}^{n}\frac{c_{i}^{\infty}}{f_{i}^{2}}\mid_{a_{\star}^{\infty},b_{\star}^{\infty}}=0$, $a_{\star}^{\infty}\neq b_{\star}^{\infty}$ then we can take $\tau_{+}^{\infty}\in (\tau_{0},b_{\star}^{\infty})$ (resp. $\tau_{-}^{\infty}\in (a_{\star}^{\infty},\tau_{0})$) such that  $\int_{\tau_{0}}^{\tau_{+}^{\infty}}h_{i}^{\epsilon=1}[c_{1}^{\infty},\ldots ,c_{n}^{\infty},D^{\infty}]=l_{i}$ (resp. $\int_{\tau_{0}}^{\tau_{-}^{\infty}}h_{i}^{\epsilon=-1}[c_{1}^{\infty},\ldots ,c_{n}^{\infty},D^{\infty}]=l_{i}$) because of the divergence of the integral towards $b_{\star}^{\infty}$ (resp. $a_{\star}^{\infty}$). Thus, by uniform convergence of $\{\sum_{i=1}^{n}\frac{c_{i}^{k}}{f_{i}^{2}}\}_{k\in \N{}}$ to $\sum_{i=1}^{n}\frac{c_{i}^{\infty}}{f_{i}^{2}}$ on compact sets we obtain the convergence of the corresponding constants $\{\tau_{\pm}^{k}\}_{k\in \N{}}$ to $\tau_{\pm}^{\infty}$ and, thus, $s_{i}^{\pm}(\hat{c}^{k},D^{k})\rightarrow s_{i}^{\pm}(\hat{c}^{\infty},D^{\infty})$. 

(iii) The remainder of the cases, except when $\frac{d}{d\tau}\sum_{i=1}^{n}\frac{c_{i}^{\infty}}{f_{i}^{2}}\mid_{\tau_{0}}=0$, $\tau_{0}=a_{\star}^{\infty}=b_{\star}^{\infty}$, are combinations of (i) and (ii).

(iv) Finally, when $\frac{d}{d\tau}\sum_{i=1}^{n}\frac{c_{i}^{\infty}}{f_{i}^{2}}\mid_{\tau_{0}}=0$, $\tau_{0}=a_{\star}^{\infty}=b_{\star}^{\infty}$, then $s_{i}^{+}(\hat{c}^{k},D^{k})\rightarrow 0$, $s_{i}^{-}(\hat{c}^{k},D^{k})\rightarrow 0$ 
(necessarily, the corresponding $K$ for $(\hat{c}^{\infty},D^{\infty})$ is $0$ and for $(\hat{c}^{k},D^{k})$ different to $0$). In fact, reasoning with $s_{i}^{+}$, the uniform convergence of $\{\sum_{i=1}^{n}\frac{c_{i}^{k}}{f_{i}^{2}}\}_{k\in \N{}}$ to $\sum_{i=1}^{n}\frac{c_{i}^{\infty}}{f_{i}^{2}}$ on compact sets implies that fixed $C_{0}>0$ and $\epsilon_{0}>0$ there exist $\eta >0$ and $k\in \N{}$ small enough and big enough respectively, such that $\int_{\tau_{0}+\eta}^{\tau_{0}+\epsilon_{0}}h_{i}^{\epsilon=1}[\hat{c}^{k},D^{k}]>C_{0}$. Thus $\int_{\tau_{0}}^{\tau_{0}+\epsilon_{0}}h_{i}^{\epsilon=1}[\hat{c}^{k},D^{k}]\rightarrow \infty$ and, thus, $s_{i}^{+}(\hat{c}^{k},D^{k})\rightarrow 0$.

For the assertion on $\mu_{i}$ define $h_{i}^{0,\epsilon}[\hat{c},D]=\epsilon\cdot \sqrt{c_{i}} \cdot f_{i}^{-2}(\tau_{0})(-D+\sum_{i=1}^{n}\frac{c_{i}}{f_{i}^{2}(\tau)})^{-1/2}$ and consider a sequence $\{(\hat{c}^{k},D^{k})\}_{k\in \N{}}$, $(\hat{c}^{k},D^{k})\rightarrow (\hat{c}^{\infty},D^{\infty})$, with $s_{1}(\hat{c}^{k},D^{k})\neq 0$ and $lim_{k\rightarrow \infty} s_{1}(\hat{c}^{k},D^{k})=0$. To simplify, assume $\epsilon=1$ and put $h_{i}^{k}\equiv h_{i}^{\epsilon=1}[\hat{c}^{k},D^{k}]$, $s_{i}^{k}\equiv s_{i}(\hat{c}^{k},D^{k})$, $h_{i}^{0,k}\equiv h_{i}^{0,\epsilon=1}[\hat{c}^{k},D^{k}]$, $s_{i}^{0,k}\equiv s_{i}^{0}(\hat{c}^{k},D^{k})$, where the functions $s_{i}^{0}$ are defined as $s_{i}$ but using $h_{i}^{0,\epsilon}$. The inequalities \begin{equation}\label{desigualdad} 
\frac{m^{k}}{M^{k}}\leq \frac{l_{i}}{\int_{\tau_{0}}^{\tau_{0}+s_{i}^{k}}h_{i}^{k}(\tau)\cdot \frac{M^{k}}{f_{i}^{-2}(\tau)}d\tau}\leq\frac{l_{i}}{\int_{\tau_{0}}^{\tau_{0}+s_{i}^{k}}h_{i}^{0,k}(\tau)d\tau}\leq \frac{l_{i}}{\int_{\tau_{0}}^{\tau_{0}+s_{i}^{k}}h_{i}^{k}(\tau)\cdot \frac{m^{k}}{f_{i}^{-2}(\tau)}d\tau}\leq \frac{M^{k}}{m^{k}}
\end{equation}
where $m^{k}$ and $M^{k}$ are, respectively, the minimum and maximum of $f_{i}^{-2}\mid [\tau_{0},\tau_{0}+s_{i}^{k}]$, imply $lim_{k\rightarrow \infty}\frac{s_{i}^{k}}{s_{i}^{0,k}}=1$ (recall $lim_{k\rightarrow \infty}s_{i}^{k}=0$). Thus $lim_{k\rightarrow \infty}\frac{s_{i}^{k}}{s_{1}^{k}}=lim_{k\rightarrow \infty}\frac{s_{i}^{0,k}}{s_{1}^{0,k}}=\frac{\sqrt{c_{1}}\cdot l_{i}\cdot f_{i}^{2}(\tau_{0})}{\sqrt{c_{i}}\cdot l_{1}\cdot f_{1}^{2}(\tau_{0})}$, as required. $\diamondsuit$

\begin{remark}\label{r1}{\rm
Clearly, reasoning as above, if $i\neq j$ then $s_{j}$ can be extended continuously to $(c_{1},\ldots ,c_{i}=0,\ldots ,c_{j}\neq 0,\ldots ,c_{n})$.
}
\end{remark}

\begin{lemma}
\label{l6}
Fix $K\in [K^{-},K^{+}]$. For each $\delta >0$, there exists $\epsilon >0$ such that if $c_{i}\leq \epsilon$ and $c_{j}\geq \delta $ for some $i$, $j$ then $\mu_{i}(\hat{c},K)<\mu_{j}(\hat{c},K)$.
\end{lemma}
{\it Proof}. Otherwise, we would find a sequence $\{\hat{c}^{m}\}$, $\hat{c}^{m}=(c_{1}^{m},\ldots ,c_{i}^{m},\ldots ,c_{j}^{m}\geq \delta ,\ldots ,c_{n}^{m})$ with $c_{i}^{m}\rightarrow 0$, $c_{1}^{m}+\cdots +c_{n}^{m}=1$ such that
\begin{equation}
\label{e26}
\mu_{i}(\hat{c}^{m},K)\geq \mu_{j}(\hat{c}^{m},K),\;\hbox{or equivalently}\; s_{i}(\hat{c}^{m},K)\leq s_{j}(\hat{c}^{m},K). 
\end{equation}
 Thus, if $\hat{c}^{\infty}=(c_{1}^{\infty},\ldots ,c_{i}^{\infty}=0,\ldots ,c_{j}^{\infty}\geq \delta ,\ldots ,c_{n}^{\infty})$ is the limit, up to a subsequence, of $\{\hat{c}^{m}\}$ then, from Remark \ref{r1}, $s_{j}(\hat{c}^{m},K)\rightarrow s_{j}(\hat{c}^{\infty},K) \in \R{}$ and we consider the following cases, which corresponds to those in the proof of Lemma \ref{l5}:

(i) If $\frac{d}{d\tau}\sum_{k=1}^{n}\frac{c_{k}^{\infty}}{f_{k}^{2}}\mid_{a_{\star}^{\infty},b_{\star}^{\infty}}\neq 0$ then $s_{i}(\hat{c}^{m},K)\rightarrow \infty >s_{j}(\hat{c}^{\infty},K)$, which contradicts (\ref{e26}).

(ii) If either $a_{\star}^{\infty}=a, b_{\star}^{\infty}=b$, or                                     $\frac{d}{d\tau}\sum_{k=1}^{n}\frac{c_{k}^{\infty}}{f_{k}^{2}}\mid_{a_{\star}^{\infty},b_{\star}^{\infty}}=0$, $a_{\star}^{\infty}\neq b_{\star}^{\infty}$, then fixed $\tau_{+}^{\infty}$, $\tau_{0}<\tau_{+}^{\infty}<b_{\star}^{\infty}$ (resp. $\tau_{-}^{\infty}$, $a_{\star}^{\infty}<\tau_{-}^{\infty}<\tau_{0}$) and taking the $\epsilon$, $D$ associated to each $(\hat{c},K)$, by using the uniform convergence in $[\tau_{0},\tau_{+}^{\infty}+\frac{b_{\star}^{\infty}-\tau_{+}^{\infty}}{2}]$ (resp. $[\tau_{-}^{\infty}-\frac{\tau_{-}^{\infty}-a_{\star}^{\infty}}{2},\tau_{0}]$) of $\{h_{i}^{\epsilon=1}[\hat{c}^{k},D^{k}]\}_{k\in \N{}}$ to $h_{i}^{\epsilon=1}[\hat{c}^{\infty},\hat{D}^{\infty}]$ (resp. $\{h_{i}^{\epsilon=-1}[\hat{c}^{k},D^{k}]\}_{k\in \N{}}$ to $h_{i}^{\epsilon=-1}[\hat{c}^{\infty},\hat{D}^{\infty}]$), we have that for $c_{i}^{k}$ small enough it satisfies $\int_{\tau_{0}}^{\tau_{+}^{\infty}}h_{i}^{\epsilon=1}[\hat{c}^{k},D^{k}]<l_{i}$ (resp. $\int_{\tau_{0}}^{\tau_{-}^{\infty}}h_{i}^{\epsilon=-1}[\hat{c}^{k},D^{k}]<l_{i}$). So, $lim_{k\rightarrow \infty}s_{i}^{+}(\hat{c}^{k},K)\geq b_{\star}^{\infty}-\tau_{0}$ (resp. $lim_{k\rightarrow \infty}s_{i}^{-}(\hat{c}^{k},K)\geq \tau_{0}-a_{\star}^{\infty}$) and thus (\ref{e26}) imply $ \int_{\tau_{0}}^{b_{\star}^{\infty}}h_{j}^{\epsilon=1}[\hat{c}^{\infty},D^{\infty}]\leq l_{j} $ (resp. $ \int_{\tau_{0}}^{a_{\star}^{\infty}}h_{j}^{\epsilon=-1}[\hat{c}^{\infty},D^{\infty}]\leq l_{j} $),   which contradicts either (\ref{e25}) if $b_{\star}^{\infty}=b$ (resp. $a_{\star}^{\infty}=a$) or  $\frac{d}{dt}\sum_{i=1}^{n}\frac{c_{i}}{f_{i}^{2}}\mid_{b_{\star}^{\infty}}=0$ if $b_{\star}^{\infty}\neq b$ (resp. $\frac{d}{dt}\sum_{i=1}^{n}\frac{c_{i}}{f_{i}^{2}}\mid_{a_{\star}^{\infty}}=0$ if $a_{\star}^{\infty}\neq a$). 

(iii) The reasoning for the remainder of the cases, except when $\frac{d}{d\tau}\sum_{k=1}^{n}\frac{c_{k}^{\infty}}{f_{k}^{2}}\mid_{\tau_{0}}=0$, $\tau_{0}=a_{\star}^{\infty}=b_{\star}^{\infty}$ are similar to (ii).

(iv) Finally, when $\frac{d}{d\tau}\sum_{k=1}^{n}\frac{c_{k}^{\infty}}{f_{k}^{2}}\mid_{\tau_{0}}=0$, $\tau_{0}=a_{\star}^{\infty}=b_{\star}^{\infty}$ (and, thus, $K=0$) then, reasoning exactly as in Lemma \ref{l5} (formula (\ref{desigualdad})), $\mu_{i}(\hat{c}^{m},0)-1+\frac{\sqrt{c_{1}^{m}}\cdot l_{i}\cdot f_{i}^{2}(\tau_{0})}{\sqrt{c_{i}^{m}}\cdot l_{1}\cdot f_{1}^{2}(\tau_{0})}\rightarrow 0$ and  $\mu_{j}(\hat{c}^{m},0)-1+\frac{\sqrt{c_{1}^{m}}\cdot l_{j}\cdot f_{j}^{2}(\tau_{0})}{\sqrt{c_{j}^{m}}\cdot l_{1}\cdot f_{1}^{2}(\tau_{0})}\rightarrow 0$ which contradicts (\ref{e26}) because $c_{i}^{\infty}=0$ and $c_{j}^{\infty}\geq \delta$. $\diamondsuit$

\begin{lemma}
\label{l7}
Consider $n$ continuous functions $\mu_{1},\ldots ,\mu_{n}:(0,1)^{n-1}\times [K^{-},K^{+}]\rightarrow \R{}$ such that $\mu_{1}\equiv 0$ and for each $\delta >0$ there exists $\epsilon >0$ satisfying $\mu_{i}(\hat{c},K)<\mu_{j}(\hat{c},K)$ whenever $c_{i}\leq \epsilon$, $c_{j}\geq \delta $.

Then there exists another $n$ continuous functions $\tilde{\mu_{1}},\ldots ,\tilde{\mu_{n}}:[-1,2]^{n-1}\times [K^{-},K^{+}]\rightarrow \R{}$, $\tilde{\mu_{1}}\equiv 0$ and constants $m_{2},\ldots ,m_{n},\bar{m_{2}},\ldots ,\bar{m_{n}}$ such that:

(i) If $\tilde{\mu_{2}}(\hat{y},K)=\cdots =\tilde{\mu_{n}}(\hat{y},K)=0$ then $\hat{y} \in (0,1)^{n-1}$ and $\mu_{2}(\hat{y},K)=\cdots =\mu_{n}(\hat{y},K)=0$.

(ii) $\tilde{\mu_{i}}(\hat{y},K)\equiv m_{i}$ for $y_{i-1}=-1$, $\tilde{\mu_{i}}(\hat{y},K)\equiv \bar{m_{i}}$ for $y_{i-1}=2$
where $m_{i}<0<\bar{m_{i}}$ for $i=2,\ldots ,n$.

(iii) There is an affine function  on $\R{}^{n-1}$ which coincide with each $\tilde{\mu_{i}} \; \; i=1,\ldots ,n$ on $\partial [-1,2]^{n-1} \times \{ K \}$ ($\partial$ denotes boundary).
\end{lemma}
{\it Proof.} From the hypothesis, there exists $0<\epsilon_{n-1}<1/2$ small enough such that 

\[
\mu_{n}(\hat{y},K)<\mu_{n-1}(\hat{y},K) \quad \hbox{if} \quad y_{n-1}(=c_{n})=\epsilon_{n-1} \quad \hbox{and} \quad \mu_{n-1}(\hat{y},K)=\cdots =\mu_{1}(\hat{y},K)(=0)
\]
\[
\mu_{n}(\hat{y},K)>\mu_{n-1}(\hat{y},K) \quad \hbox{if} \quad y_{n-1}=1-\epsilon_{n-1} \quad \hbox{and} \quad \mu_{n-1}(\hat{y},K)=\cdots =\mu_{1}(\hat{y},K)(=0) 
\]

Moreover, fixed $t\in \{1,\ldots ,n-2\}$ and constants $\epsilon_{t+1},\ldots ,\epsilon_{n-1}\in (0,1/2)$, there exists $\epsilon_{t}\in (0,1/2)$ such that

\[
\begin{array}{lr}
\mu_{t+1}(\hat{y},K)<\mu_{t}(\hat{y},K) & \hbox{if} \quad y_{t}=\frac{c_{t+1}}{c_{1}+\cdots +c_{t+1}}=\epsilon_{t}, (y_{t+1},\ldots ,y_{n-1})\in \prod_{j=t+1}^{n-1}[\epsilon_{j},1-\epsilon_{j}] \\                                                                                 
                                    & \hbox{and} \quad \mu_{t}(\hat{y},K)=\cdots =\mu_{1}(\hat{y},K)(=0)
\end{array}
\]
\begin{equation}\label{penultima}
\end{equation}
\[
\begin{array}{lr}
\mu_{t+1}(\hat{y},K)>\mu_{t}(\hat{y},K) & \hbox{if} \quad y_{t}=\frac{c_{t+1}}{c_{1}+\cdots +c_{t+1}}=1-\epsilon_{t}, (y_{t+1},\ldots ,y_{n-1})\in \prod_{j=t+1}^{n-1}[\epsilon_{j},1- \epsilon_{j}] \\   
                                    & \hbox{and} \quad \mu_{t}(\hat{y},K)=\cdots =\mu_{1}(\hat{y},K)(=0)
\end{array}
\]

This allows us to obtain, by induction, an $\epsilon^{\star}>0$ such that the restriction of the functions $\mu_{1},\ldots ,\mu_{n}$ to $[\epsilon^{\star},1-\epsilon^{\star}]^{n-1}\times[K^{-},K^{+}]$ satisfy:

\begin{equation}\label{ultima}\begin{array}{c}
\mu_{i}(\hat{y},K)<\mu_{i-1}(\hat{y},K) \quad \hbox{if} \quad y_{i-1}=\epsilon^{\star} \quad \hbox{and} \quad \mu_{i-1}(\hat{y},K)=\cdots =\mu_{1}(\hat{y},K) \\ \mu_{i}(\hat{y},K)>\mu_{i-1}(\hat{y},K) \quad \hbox{if} \quad y_{i-1}=1-\epsilon^{\star} \quad \hbox{and} \quad \mu_{i-1}(\hat{y},K)=\cdots =\mu_{1}(\hat{y},K)
\end{array}
\end{equation}

Extend continuously these functions to $[0,1]^{n-1}\times [K^{-},K^{+}]$ by the following procedure of $n-1$ steps. The $(i-1)-th$ step, $i=2,\ldots ,n$ is carried out as follows. 

(A) Consider the functions $\mu_{1},\ldots ,\mu_{i-1},\mu_{i+1},\ldots ,\mu_{n}$ and extend them trivially in the variable $y_{i-1}$ (that is, independently of this variable) until $y_{i-1}=0$ and until $y_{i-1}=1$.

(B) Consider the function $\mu_{i}$ and extend it linearly (in the variable $y_{i-1}$) until $y_{i-1}=0$  such that $\mu_{i}\equiv v_{i}<0$ ($v_{i}$ arbitrary) for $y_{i-1}=0$. Analogously, extend the function $\mu_{i}$ linearly (in the variable $y_{i-1}$) until $y_{i-1}=1$ in such a way that $\mu_{i}\equiv \bar{v_{i}}>0$ for $y_{i-1}=1$.

Finally, to obtain {\it (iii)} too, extend continuously these functions to $[-1,2]^{n-1}\times [K^{-},K^{+}]$ by the following process.
First extend the function $\mu_{1}$ as $0$ in all the variables (thus $m_{1}=0$). Then consider the function $\mu_{i}$ for $i=2,\ldots ,n$ and extend it linearly in the variable $y_{i-1}$ until $y_{i-1}=-1$ and until $y_{i-1}=2$ respectively such that $\mu_{i}\equiv 2v_{i}-\bar{v}_{i}$ for $y_{i-1}=-1$ and $\mu_{i}\equiv 2\bar{v}_{i}-v_{i}$ for $y_{i-1}=2$ (thus $m_{i}=2v_{i}-\bar{v}_{i}$, $\bar{m}_{i}=2\bar{v}_{i}-v_{i}$). At last, extend the function $\mu_{i}$ linearly in the variables $y_{j}$ for all $j\neq i-1$ in such a way that $\mu_{i}(\hat{y},K)=m_{i}+\frac{\bar{m}_{i}-m_{i}}{3}(y_{i-1}+1)$ for $y_{j}\in \{-1,2\}$. $\diamondsuit$
  

\section{Proof of Theorem \ref{t1}}

\begin{theorem}
\label{t4}
A spacetime $(I\times F_{1}\times \cdots \times F_{n},g)$ with weakly convex fibers $(F_{1},g_{1}),\cdots ,(F_{n},g_{n})$ and such that

\begin{equation}\label{e28} \begin{array}{ll} \int_{c}^{b}f_{i}^{-2}(\frac{1}{f_{1}^{2}}+\cdots +\frac{1}{f_{n}^{2}})^{-1/2}=\infty  & \int_{a}^{c}f_{i}^{-2}(\frac{1}{f_{1}^{2}}+\cdots +\frac{1}{f_{n}^{2}})^{-1/2}=\infty  \end{array} \end{equation}

for all $i$ and for $c\in (a,b)$ is geodesically connected.
\end{theorem}

As equalities (\ref{e28}) are weaker than (\ref{e24}) now the analogous to Lemma \ref{l4} is weaker too.

\begin{lemma}
\label{l9}
Property (\ref{e28}) imply (\ref{e25}) for $D\geq 0$ and for all $c_{1},\ldots ,c_{n}\geq 0$, $c_{1}+\cdots +c_{n}=1$ where it makes sense. Moreover, if one of the integrals (\ref{e25}) is finite for one $D<0$ then it is finite for all $D<0$ and, at any case, the values of the integrals (in $]0,+\infty]$) depends  continuously of $\hat{c},D$.
\end{lemma}

 The main difference between the proof of Theorem \ref{t4} and Theorem \ref{t3} is that, now, Lemma \ref{l4} does not hold for $D<0$ and, thus perhaps $b_{\star}=b$ but $\int_{\tau_{0}}^{b_{\star}}h_{i}^{\epsilon=1}<l_{i}$ for all $i$. In this case Convention \ref{c1} was not appliable because, otherwise, the corresponding function $\tau(t_{i}(r))$ defined from (\ref{e4}) and (\ref{e5}) touches an extreme of $I=(a,b)$ before reaching $z'$ and, thus, it does not correspond with the projection of a (true) geodesic in $(a,b)\times F_{1}\times \cdots \times F_{n}$ joining $z$ and $z'$. Nevertheless, we will admit this possibility in order to extend the arguments on continuity in Theorem \ref{t3}, and $\tau(t_{i}(r))$ will be regarded as the projection of a {\it fake} geodesic in $[a,b]\times F_{1}\times \cdots \times F_{n}$. Then, an additional effort will be necessary to ensure that, among the obtained {\it generalized} geodesics, enough true geodesics are yielded. Even more, we will have to admit fake geodesics even when an extreme of $I$ is infinite so, an auxiliar diffeomorphism will be used to normalize. 

{\it Proof of Theorem \ref{t4}.} Fix a diffeomorphism $\varphi :(a,b)\rightarrow (-1,1)$ $\varphi'>0$ and redefine the functions $s_{i}^{\pm}$ by composing all the values in $I$ with $\varphi$ and making use of Convention \ref{c1} for each $i$, even if $a_{\star}=a$ or $b_{\star}=b$, that is: $s_{i}^{+}(\hat{c},K)$ is equal to $\varphi(t_{i})-\varphi(\tau_{0})$, $2\varphi(b^{\star})-\varphi(t_{i})-\varphi(\tau_{0})$, etc. depending on the different cases of the integral of $h_{i}^{+}$ (analogously with $s_{i}^{-}$). The continuity of these functions $s_{i}^{\pm}$ (thus, the conclusion on $\mu_{i}$) can be proven as in Lemma \ref{l5}, but taking into account that now, in the case (ii) when $b_{\star}^{\infty}=b$ perhaps $\int_{\tau_{0}}^{b_{\star}^{\infty}}h_{i}^{\epsilon=1}[\hat{c}^{\infty},D^{\infty}]<\infty$ and, so, the continuity assertion in Lemma \ref{l9} must be used. On the other hand, the boundary conditions in Lemma \ref{l6} must be regarded as follows. Put  $\bar{K}^{+}=min_{i}\{\frac{1}{f_{i}^{2}(\tau_{0})}\}$, $\bar{K}^{-}=-\bar{K}^{+}$; note that when $K\in [\bar{K}^{-},\bar{K}^{+}]$  necessarily $D\geq 0$ and, in this particular case, Lemma \ref{l9} can be claimed. So, Lemma \ref{l6} still holds if $[K^{-},K^{+}]$ is replaced by $[\bar{K}^{-},\bar{K}^{+}]$. This modification makes necessary the following longer version of  Lemma \ref{lA} in the proof of Theorem \ref{t3}.  
                      
\begin{lemma} \label{lB}
The set of points $(\hat{y},K)\in (0,1)^{n-1}\times [K^{-},K^{+}]$ verifying $\mu_{2}(\hat{y},K)=\cdots =\mu_{n}(\hat{y},K)=0$ admits a connected component      ${\cal C}$ verifying one of the following possibilities:

(i) ${\cal C}\cap \{(0,1)^{n-1}\times \{K^{+}\}\}\neq \emptyset$, ${\cal C}\cap \{(0,1)^{n-1}\times \{K^{-}\}\} \neq \emptyset$.

(ii) ${\cal C}\cap \{(0,1)^{n-1}\times \{K^{+}\}\}\neq \emptyset$, ${\cal C}\cap \{\partial[0,1]^{n-1}\times [K^{-},\bar{K^{-}}]\}\neq \emptyset$.

(iii) ${\cal C}\cap \{\partial[0,1]^{n-1}\times [\bar{K^{+}},K^{+}]\}\neq \emptyset$, ${\cal C}\cap \{(0,1)^{n-1}\times \{K^{-}\}\} \neq \emptyset$.

(iv) ${\cal C}\cap \{\partial[0,1]^{n-1}\times [\bar{K^{+}},K^{+}]\}\neq \emptyset$, ${\cal C}\cap \{\partial[0,1]^{n-1}\times [K^{-},\bar{K^{-}}]\}\neq \emptyset$.
\end{lemma}

{\it Proof of Lemma \ref{lB}.} By an analogous reasoning to Lemma \ref{l7} and for $\epsilon^*$ small enough the restriction of the functions $\mu_{1},\ldots ,\mu_{n}$ to $[\epsilon^*,1-\epsilon^*]^{n-1}\times [\bar{K^{-}},\bar{K^{+}}]$ satisfy:
\begin{equation}
\label{e28.5}
\begin{array}{c}
 \mu_{i}(\hat{y},K)<\mu_{i-1}(\hat{y},K) \quad \hbox{if} \quad y_{i-1}=\epsilon^* \quad \hbox{and} \quad \mu_{i-1}(\hat{y},K)=\cdots =\mu_{1}(\hat{y},K)(=0) \\
\mu_{i}(\hat{y},K)>\mu_{i-1}(\hat{y},K) \quad \hbox{if} \quad y_{i-1}=1-\epsilon^* \quad \hbox{and} \quad \mu_{i-1}(\hat{y},K)=\cdots =\mu_{1}(\hat{y},K)(=0) 
\end{array}
\end{equation}
for $i=2,\ldots ,n$. Now, consider a homeomorphism $h:[\epsilon^*,1-\epsilon^*]^{n-1}\times [K^{-},K^{+}]\rightarrow [\epsilon^*,1-\epsilon^*]^{n-1}\times [-1,1]$ such that $\quad h(\partial [\epsilon^*,1-\epsilon^*]^{n-1}\times [\bar{K}^{-},\bar{K}^{+}])=\partial [\epsilon^*,1-\epsilon^*]^{n-1}\times [-1,1]$, $\quad h(\partial [\epsilon^*,1-\epsilon^*]^{n-1}\times [\bar{K}^{+},K^{+}]\cup [\epsilon^*,1-\epsilon^*]^{n-1}\times \{K^{+}\})=[\epsilon^*,1-\epsilon^*]^{n-1}\times \{1\}$, $\quad h(\partial [\epsilon^*,1-\epsilon^*]^{n-1}\times [K^{-},\bar{K}^{-}]\cup [\epsilon^*,1-\epsilon^*]^{n-1}\times \{K^{-}\})=[\epsilon^*,1-\epsilon^*]^{n-1}\times \{-1\}$, clearly, functions $\mu_{i}\circ h^{-1}$ satisfy (\ref{e28.5}) but in $[\epsilon^*,1-\epsilon^*]^{n-1}\times [-1,1]$. Following the reasoning in Lemma \ref{l7}, from the functions $\mu_{i}\circ h^{-1}$ we obtain functions $\widetilde{\mu_{i}\circ h^{-1}}$ where the arguments of degree of solutions can be applied as in Lemma \ref{lA}. By applying standard arguments (\cite[Ch.I Theor. 9.1]{Wh})      to the corresponding connected components obtained from $\{\epsilon^*_{n}\}_{n\in \N{}}$, $\epsilon^*_{n}\rightarrow 0$, we obtain a connected component which composed with $h^{-1}$ clearly satisfies one of the conditions {\it (i)-(iv)} of Lemma \ref{lB}. $\diamondsuit$

Now, each point of ${\cal C}$ represents either (i) a true geodesic reaching the line $L[x_{1}',\ldots ,x_{n}']$ or (ii) a fake geodesic in $[a,b]\times F_{1}\times \cdots \times F_{n}$ that touches $\tau=b$ or $\tau=a$ before reaching the line $L[x_{1}',\ldots ,x_{n}']$. Recall that such a fake geodesic $\gamma(t)=(\tau(t),\gamma_{1}(t),\ldots ,\gamma_{n}(t))$ can touch either an extreme or both extremes of $[a,b]$. Then, $\gamma$ will be said a fake geodesic at $b$ (resp. $a$) if $b$ (resp. $a$) is the first extreme touched. In this case, if $t_{0}$ is the first point such that $\tau(t_{0})=b$ (resp. $a$) $\gamma(t_{0})$ will be called the {\it escape point} of the fake geodesic.

The timelike component $\tau(t)$ of a (true or fake) geodesic depends on $\hat{y}$ and $K$ continuously. Since ${\cal C}$ is connected, if we suppose that every point of ${\cal C}$ represents a true geodesic then the subset $J$ of the line $L[x_{1}',\ldots,x_{n}']$ which can be reached by true geodesics in ${\cal C}$ is an interval. Moreover, the subset $J$ is also an interval even if there exists some fake geodesic by Lemma \ref{l10}.

Now, we consider two cases:

(i) If ${\cal C}$ satisfies {\it (i)} in Lemma \ref{lB}, then by the choice of $K^{+}$,$K^{-}$ either $\tau_{0}'\in J$ and the proof is over or $z'\in J^{+}(z)\cup J^{-}(z)$ and, in this case, Theorem \ref{t2} can be claimed.

(ii) If ${\cal C}$ belongs to the case {\it (ii)} (resp. {\it (iii)};{\it (iv)}) then $J=[a',b']$ with possibly $a'=a$ or $b'=b$ (and in this case $[a',b']$ must be assumed open in the corresponding extreme). Moreover, we have that $a'=a$ (resp. $b'=b$; $a'= a$,$b'=b$). In fact, let $(c^{\infty},D^{\infty})$ be a point in ${\cal C}\cap(\partial [0,1]^{n-1}\times ]K^{-},K^{+}[)$ and $\{(c^{k},D^{k})\}_{k\in \N{}}$ a sequence in ${\cal C}$ convergent to this point. Recall that as the sequence lies in ${\cal C}$ the conclusion of Lemma \ref{l6} cannot hold. This implies that the sequence must satisfy (ii) of this Lemma with $a_{\star}^{\infty}=a$, $b_{\star}^{\infty}=b$. If $\int_{\tau_{0}}^{a}(c
_{i}^{\infty})^{-1/2}\cdot h_{i}^{\epsilon=-1}[\hat{c}^{\infty},D^{\infty}]<\infty$ (resp. $\int_{\tau_{0}}^{b}(c_{i}^{\infty})^{-1/2}\cdot h_{i}^{\epsilon=1}[\hat{c}^{\infty},D^{\infty}]<\infty$; both inequalities hold) for some $i$ such that $c_{i}^{\infty}=0$ then for $k$ big enough each $(\hat{c}^{k},D^{k})$ corresponds with a fake geodesic and Lemma \ref{l10} can be claimed (recall that $h_i^{\epsilon}$ contains a factor $(c_i^\infty)^{1/2}$ which maybe $0$; so, we have multiplied by $(c_i^\infty)^{-1/2}$ to remove this factor). Otherwise, the result can be obviously obtained. Thus, either $\tau_{0}'\in J$ and the proof is over or, from the choice of $K^{+}$, $K^{-}$, we have $z'\in J^{+}(z)\cup J^{-}(z)$ and Theorem \ref{t2} can be claimed. $\diamondsuit$ 

\begin{lemma}
\label{l10}
If there exists some point $p\in {\it {\cal C}}$ which represents a fake geodesic at $b$ (resp. $a$) then the set of points $J$ of $L[x_{1}',\ldots ,x_{n}']$ reached by true geodesics is an interval closed in $(a,b)$ with an extreme equal to $b$ (resp. $a$).
\end{lemma}
{\it Proof.} First note that for $K=0$ necessarily $D>0$, thus Lemma \ref{l9} is appliable and all these points in ${\cal C}$ represents true geodesics and so $J\neq \emptyset$. If $J$ were not such an interval, then there would exist $\tau_{0}''\in (a,b)-J$ such that $\tau_{0}''-\delta_{1}\in J$, $\delta_{1}>0$. Define the function which  maps every point $(\hat{c},K)\in {\cal C}$ into either the escape point of $\gamma(\hat{c},K)$ if it is a fake geodesic or the point where $L[x_{1}',\ldots ,x_{n}']$ is reached by $\gamma$ otherwise. The restriction of this function to the set of all the true geodesics of ${\cal C}$ is continuous because of Lemma \ref{l5} and so is its restriction to the (closed) set of all the fake geodesics in ${\cal C}$, because of Lemma \ref{l9}. If $(c^{\infty},D^{\infty})$ represents a fake geodesic in ${\cal C}$ and $\{(\hat{c}^{k},D^{k})\}$ a sequence of true geodesics in ${\cal C}$ convergent to it, necessarily $D^{\infty}<0$ and, from Lemma \ref{l9}, $\int_{\tau_{0}}^{b}h_{i}^{k}\rightarrow \int_{\tau_{0}}^{b}h_{i}^{\infty}\leq l_{i}$ for all $i$. As $l_{i}=\int_{\tau_{0}}^{t_{i}^{k}}h_{i}^{k}\leq \int_{\tau_{0}}^{b}h_{i}^{k}$ for some $t_{i}^{k}<b$, necessarily $\{t_{i}^{k}\}\rightarrow b$ for all $i$ and the continuity is obtained. Thus the image of the defined function on ${\cal C}$                                                                     is a connected subset of $[a,b]\times F_{1}\times \cdots \times F_{n}$. But this is a contradiction because  $[a,\tau_{0}''[\times F_{1}\times \cdots \times F_{n}$ and $]\tau_{0}'',b]\times F_{1}\times \cdots \times F_{n}$                           contain points of the image. Moreover, from the continuity of the previous function and the compactness of ${\cal C}$, $J$ is closed in $(a,b)$. $\diamondsuit$

Finally, Theorem \ref{t1} is a consequence of Theorem \ref{t4} and the  proposition below the following lemma.

\begin{lemma}\label{lemma} (1) Let $\phi^{k}$, $\phi$ be positive functions on the interval $[d,b[$ with $\{\phi^{k}\}\rightarrow \phi$ uniformly on compact subsets. If 

\begin{equation}\label{tres} lim sup\int_{d}^{b}\phi^{k}\leq \int_{d}^{b}\phi\in ]0,\infty]\end{equation}
then $lim\int_{d}^{b}\phi^{k}$ exists and is equal to $\int_{d}^{b}\phi$.

(2) If $\hat{c}^{k}, \hat{c}^{\infty}\in \overline{\bigtriangleup_{n}}=\{(c_{1},\ldots ,c_{n})\in [0,1]: \sum_{i}c_{i}=1\}$ and $\{\hat{c}^{k}\}\rightarrow\hat{c}^{\infty}$ then $lim \int_{d}^{b}f_{i}^{-2}(\frac{c_{1}^{k}}{f_{1}^{2}}+\cdots +\frac{c_{n}^{k}}{f_{n}^{2}})^{-1/2}=
\int_{d}^{b}f_{i}^{-2}(\frac{c_{1}^{\infty}}{f_{1}^{2}}+\cdots +\frac{c_{n}^{\infty}}{f_{n}^{2}})^{-1/2}$.
\end{lemma}

{\it Proof.} (1) In the finite case note that for $\eta$ small $\int_{d}^{b}\phi-\eta=\int_{d}^{b-\delta}\phi=lim\int_{d}^{b-\delta}\phi^{k}\leq lim inf\int_{d}^{b}\phi^{k}$ for some $\delta>0$. The infinite case is analogous.

(2) It is sufficient to check that inequality (\ref{tres}) holds. For any $\lambda\in ]0,1[$ put $d_{i}^{\infty}=\lambda^{2}c_{i}^{\infty}$ and take $k_{0}$ such that if $k\geq k_{0}$ then $c_{i}^{k}\geq d_{i}^{\infty}$ for all $i$. Then  
\[
 \int_{d}^{b}f_{i}^{-2}(\frac{c_{1}^{k}}{f_{1}^{2}}+\cdots +\frac{c_{n}^{k}}{f_{n}^{2}})^{-1/2}\leq
\int_{d}^{b}f_{i}^{-2}(\frac{d_{1}^{\infty}}{f_{1}^{2}}+\cdots +\frac{d_{n}^{\infty}}{f_{n}^{2}})^{-1/2}=\frac{1}{\lambda}\cdot\int_{d}^{b}f_{i}^{-2}(\frac{c_{1}^{\infty}}{f_{1}^{2}}+\cdots +\frac{c_{n}^{\infty}}{f_{n}^{2}})^{-1/2}.
\]
Thus, the result follow making $\lambda\rightarrow 1$. $\diamondsuit$

\begin{proposition}\label{p1} 
Any point of the spacetime can be joined with any line by means of both, a future directed and a past directed causal curve if and only if equalities (\ref{e28}) hold.
\end{proposition}

{\it Proof.} Consider the point $(\tau,x_{1},\ldots ,x_{n})$ and the line  $L[x'_{1},\ldots ,x'_{n}]$ with $x_{i}'\neq x_{i}$ for all $i$. By hypothesis, there exists $\tau'>\tau$ such that $(\tau',x'_{1},\ldots ,x'_{n}) \in J^{+}(\tau,x_{1},\ldots ,x_{n})$. By Theorem \ref{t2} we can find $c_{1},\ldots ,c_{n}\geq 0$, $c_{1}+\cdots +c_{n}=1$ such that

\[ \sqrt{c_{i}}\int_{\tau}^{\tau'}f_{i}^{-2}(\frac{c_{1}}{f_{1}^{2}}+\cdots +\frac{c_{n}}{f_{n}^{2}})^{-1/2}\geq l_{i} \] 
for all i.

Repeating the procedure, making $\tau^{k}=b-1/k$, we obtain $\tau'(\tau^{k})\rightarrow b$, $\tau^{k}<\tau'(\tau^{k})$ and constants $(c_{1}^{k}(\tau),\ldots ,c_{n}^{k}(\tau))\rightarrow (c_{1}^{\infty},\ldots ,c_{n}^{\infty})\in \overline{\bigtriangleup_{n}}$ (up to a subsequence). Now if              

\begin{equation}\label{uno} \sqrt{c_{i}^{\infty}}\int_{c}^{b}f_{i}^{-2}(\frac{c_{1}^{\infty}}{f_{1}^{2}}+\cdots +\frac{c_{n}^{\infty}}{f_{n}^{2}})^{-1/2}=L<\infty
\end{equation}
then for some $\delta>0$

\begin{equation}\label{dos} \sqrt{c_{i}^{\infty}}\int_{b-\delta}^{b}f_{i}^{-2}(\frac{c_{1}^{\infty}}{f_{1}^{2}}+\cdots +\frac{c_{n}^{\infty}}{f_{n}^{2}})^{-1/2}=l_{i}.
\end{equation}
Thus, taking $k>1/\delta$, Lemma \ref{lemma} (2)  with $d=b-\delta$ and formulae (\ref{uno}), (\ref{dos}) yield a contradiction. $\diamondsuit$

\section{A more general version and applications.} \label{s6}

Analizing the proof of the Theorem in previous section,  we note that the behaviour of the functions 
\begin{equation} \label{fh}
f^{\hat c}(\tau) = \sum_{i=1}^{n}\frac{c_{i}}{f_{i}^{2}(\tau)}
\end{equation}
 at the extremes of the interval $I$, is the essential point for the geodesic connectedness. In this sense, we give the following:

\begin{definition}\label{d1} Let $f:(a,b)\rightarrow \R{}$ be a smooth function and let $m_{b}=lim inf_{\tau\rightarrow b}f(\tau)$ (resp. $m_{a}=lim inf_{\tau\rightarrow a}f(\tau)$). The extreme $b$ (resp. $a$) is a (strict) relative minimum of $f$ if there exists $\epsilon>0$ such that if $0<\delta_{1}<\epsilon$, then $f(b-\delta_{1})>m_{b}$ (resp. $f(a+\delta_{1})>m_{a}$).  
\end{definition}

Now, we are ready to display the hipothesis which will allow to obtain geodesic connectedness.

{\bf Condition} (*). {\it For every $c_{1},\ldots ,c_{n}\geq 0$, $c_{1}+\cdots +c_{n}=1$ such that the function $f^{\hat c}$ in (\ref{fh}) reaches in $b$ (resp. in $a$) a relative minimum it verifies $\int_{c}^{b}f_{i}^{-2}(\frac{c_{1}}{f_{1}^{2}}+\cdots +\frac{c_{n}}{f_{n}^{2}}-m_{b})^{-1/2}=\infty$ (resp. $\int_{a}^{c}f_{i}^{-2}(\frac{c_{1}}{f_{1}^{2}}+\cdots +\frac{c_{n}}{f_{n}^{2}}-m_{a})^{-1/2}=\infty$) for all $i$ and for some $c\in (a,b)$. }

[Recall that if one of the previous integrals is finite for one $D<m_{a}$ (resp. $D<m_{b}$) then it is finite for all $D<m_{a}$ (resp. $D<m_{b}$) and, at any case, the values of the integrals (in $]0,+\infty]$) depends continuously on $\hat{c}$, $D$].

In fact, Condition (*) makes sure that for every $c_{1},\ldots ,c_{n}\geq 0$, $c_{1}+\cdots +c_{n}=1$,  the integral of $h_{i}^{\epsilon}[\hat{c},D]$ around every relative minimum (ordinary or in the sense of Definition \ref{d1} )  is arbitrarily big when $D$ is close to the value of this minimum, which will be sufficient as we will check.

The main differences between the proofs for this new case and Theorem \ref{t4} are the following.

Clearly, Condition (*) can replace Lemma \ref{l9} (or Lemma \ref{l4}) in the proof of Lemma \ref{l5} and Remark \ref{r1}, but the modifications for Lemmas \ref{l6}, \ref{l7} are not so simple, and the  following construction will be needed. Replace the set $\{(\hat{c},K): \hat{c}\in\overline{\bigtriangleup}_{n}, K\in [\bar{K}^{-},\bar{K}^{+}]\}$ by the set $G_{N,\eta}=\{(\hat{c},K): \hat{c}\in\overline{\bigtriangleup}_{n}, K\in [K_{N}(\hat{c})-\eta, K_{N}(\hat{c})+\eta]\}, N\in \N{}, \eta>0$ where $K_{N}(\hat{c})=sign(m_{r}(\hat{c})-m_{l}(\hat{c}))\cdot (f^{\hat{c}}(\tau_{0})-D_{N}(\hat{c}))$ being $D_{N}(\hat{c})=N\cdot e(\hat{c})\cdot |m_{l}(\hat{c})-m_{r}(\hat{c})|  \cdot (m(\hat{c})-f^{\hat{c}}(\tau_{0}))+f^{\hat{c}}(\tau_{0})$ with $e:\bigtriangleup_{n}\rightarrow \R{}$,

\begin{equation} e(\hat{c})=\left\{ \begin{array}{lr} 1 & \hbox{if} \mid m_{l}(\hat{c})-m_{r}(\hat{c}) \mid < 1/N \\ \frac{1}{N\cdot \mid m_{l}(\hat{c})-m_{r}(\hat{c}) \mid} & \hbox{if} \mid m_{l}(\hat{c})-m_{r}(\hat{c}) \mid \geq 1/N \end{array} \right.,
\end{equation}
$m_{l}(\hat{c})$ the minimum value of $f^{\hat{c}}$ to the left of $\tau_{0}$, $m_{r}(\hat{c})$ the minimum value of $f^{\hat{c}}$ to the right of $\tau_{0}$, $m(\hat{c})=min\{m_{l}(\hat{c}),m_{r}(\hat{c})\}$. $N$, $\eta$ will be chosen large enough and small enough, respectively, in order to obtain, by induction, $\epsilon^{\star}>0$ such that the relations analogous to (\ref{ultima}) holds (the inequalities of these relations are now written as $\mu_{i}(\hat{c},K_{N}(\hat{c})+\eta')<\mu_{j}(\hat{c},K_{N}(\hat{c})+\eta')$). In fact, $N$ and $\eta$ can be taken as the maximum and the minimum, respectively, of the $n-1$ ones obtained as follows in each step of the analog to the inductive proceadure of Lemma \ref{l7}. In the ($n-t$)-th step, given the relations analogous to (\ref{penultima}), $\eta$ and $N$ are  such that whenever $\hat{c}\in\overline{\bigtriangleup}_{n}$, $(y_{t+1},\ldots ,y_{n-1})\in \prod_{j=t+1}^{n-1}[\epsilon_{j},1-\epsilon_{j}]$, $N'>N$, $\eta'\in[-\eta,\eta]$ satisfy  either $a_{\star}=a$ or $b_{\star}=b$ then $\int_{\tau_{0}}h_{j}^{\epsilon}[\hat{c},K_{N'}(\hat{c})+\eta']>l_{j}$ for some $j\leq t+1$, where the upper limit of the integral is $b_{\star}$ if $\epsilon(=sign(K_{N'}(\hat{c})+\eta'))=1$ and $a_{\star}$ if $\epsilon=-1$. Then,  a suitable version of Lemma \ref{lB} is obtained, say: the set of points $(\hat{y},K)\in (0,1)^{n-1}\times [K^{-},K^{+}]$ verifying $\mu_{2}(\hat{y},K)=\cdots =\mu_{n}(\hat{y},K)=0$ admits a connected component ${\cal C}$ verifying one of the possibilities in Lemma \ref{lB} but replacing $\{\partial[0,1]^{n-1}\times [K^{-},\bar{K^{-}}]\}$ and $\{\partial[0,1]^{n-1}\times [\bar{K^{+}},K^{+}]\}$ by  $\{(\hat{c},K):\hat{c}\in \partial[0,1]^{n-1}, K\in [K^{-}, K_{N}(\hat{c})-\eta]\}$ and $\{(\hat{c},K):\hat{c}\in \partial[0,1]^{n-1}, K\in [ K_{N}(\hat{c})+\eta, K^{+}]\}$ respectively. From here the proof goes on as before taking into account that Condition (*) is used instead of Lemma \ref{l9} in the proof of Lemma \ref{l10} and the following subtlety. For the proof of, say, $a'=a$ we knew that $\epsilon=-1$ because the value $K^{\infty}$ of $K$ corresponding to $(\hat{c}^{\infty},D^{\infty})$ was negative. Now, it is possible for $K^{\infty}$ to be positive, $0<K^{\infty}<K_{N}(\hat{c}^{\infty})$ but, in this case the definition of $K_{N}$ implies $b_{\star}^{\infty}<b$. 

Thus, we obtain:

\begin{theorem}\label{t5} A multiwarped spacetime $(I\times F_{1}\times \cdots \times F_{n},g)$ with weakly convex fibers $(F_{1},g_{1}),\ldots ,(F_{n},g_{n})$ verifying Condition (*) above is geodesically connected.
\end{theorem}

Next, we will apply this result to some well-known spacetimes.                  
First, note that the (geodesically connected) strip of Lorentz-Minkowski spacetime $\L{}^{n+1}$ commented in the Introduction satisfies Condition (*), of course. Less trivially, consider the metric on $(r_{-} ,r_{+})\times \R{} \times S^{2}$
\[
       - \left( 1 -\frac{2m}{r} + \frac{e^{2}}{r^{2}} \right) ds^{2}  + 
\left( 1- \frac{2m}{r} +\frac{e^{2}}{r^{2}} \right)^{-1}  dr^{2}  + 
r^{2} (d\theta ^{2}  + \sin ^{2} \theta d\phi ^{2})
\]
where $e, m $ are constant ($m$ positive), $e^{2} < m^{2}, \quad r\in (r_{-} ,r_{+}),  \quad r_{-} =m-(m^{2} -e^{2} )^{1/2}   , 
\quad r_{+} =m+(m^{2} -e^{2})^{1/2}   , s\in \R{}$ and $\theta, \phi$ are usual coordinates on $S^{2}$  (note  that  the  new variable $\tau\in (0,\tau_{+})$ obtained integrating 
\[
  d\tau  = \left(-1 +\frac{2m}{r} - \frac{e^{2}}{r^{2}} \right)^{-1/2} dr  
\]
yields the interval $I$ and the warping functions  of  our  definition  of  multiwarped 
spacetime; on the other hand,  the fibers are, of course, weakly convex).

When  $e\neq 0$  this  metric  is the intermediate zone  of  usual  Reissner-Nordström spacetime and, in this case,  $\frac{1}{f_{1}^{2}(\tau)}$ diverges on the extremes of the interval and  $\frac{1}{f_{2}^{2}(\tau)}$       has relative minima which satisfy Condition (*). Therefore, Condition (*) always hold and it is geodesically connected.

When $e=0$ this metric represents the inner piece of Schwarzschild spacetime, and functions $f_{1}(\tau)$, $f_{2}(\tau)$ behaves at $0$ as $\tau^{-1/3}$, $\tau^{2/3}$  and at $\tau_{+}$ as $\tau_{+}-\tau$, $(\tau_{+}-(\tau_{+}-\tau)^{2})^{2}$, thus Condition (*) does not hold for $f_{1}$ at $0$. Nevertheless, its geodesic connectedness is a consequence of the special characteristics of its fiber $F_{2}$. To simplify things, the following result will be stated just for $n=2$ fibers, so $c_{1}\equiv c$, $c_{2}\equiv 1-c$.

\begin{lemma}\label{lemma'} Let $(I\times F_{1}\times F_{2},g)$ be a multiwarped spacetime with weakly convex fibers $(F_{1},g_{1})$, $(F_{2},g_{2})$ and assume Condition (*) holds for all $c\in [0,1]$ but when $c=1$, $i=1$. Then $(I\times F_{1}\times F_{2},g)$ is geodesically connected if the second fiber satisfies the following condition: for all $x_{2}, x'_{2}\in F_{2}$ there exist a sequence of geodesics of the fiber $\{\overline{\gamma}_{2}^{m}\}_{m\in \N{}}$ joining them with diverging lenghts $L_{m}$ (this condition happens if each pair $x_{2},x_{2}'\in F_{2}$ is crossed by a closed geodesic).
\end{lemma}

{\it Proof.} Given $z=(t, x_{1},x_{2})\in I\times F_{1}\times F_{2}$ and the line $L[x_{1}',x_{2}']$ fix a minimizing geodesic $\overline{\gamma}_{1}$ joining $x_{1}$ and $x'_{1}$ and consider $s_{1}(c,K)$, $c\in ]0,1]$. For each $m\in \N{}$, the corresponding $\overline{\gamma}_{2}^{m}$ allows to define the function $s_{2}^{m}(c,K)$, $c\in [0,1[$. Assuming that $\{L_{m}\}_{m\in \N{}}$ is increasing then $\{s_{2}^{m}(c,K)\}_{m\in \N{}}$ is increasing too. 

Fixed $\rho>0$, there exist $m$ such that $0 \equiv \mu_{1}(1-\rho,K_{N}(1-\rho))<\mu_{2}^{m}(1-\rho,K_{N}(1-\rho))$ for all $N$ big enough. Moreover, as in the proof of Theorem \ref{t5}, there exist $N$,$\eta$ large enough and small enough, respectively, such that  $0\equiv \mu_{1}(\epsilon^*,K_{N}(\epsilon^*)+\eta')>\mu_{2}^{m}(\epsilon^*,K_{N}(\epsilon^*)+\eta')$ holds for any $\eta'\in [-\eta,\eta]$,  and  $\epsilon^*$ small enough. Then,  we obtain a connected component ${\cal C}_{\rho}$ satisfying one of the corresponding four possibilities in Lemma \ref{lB} (replacing too $(0,1)\equiv (0,1)^{n-1}$ by $(0,1-\rho)$).

Thus we can connect $z$ with an interval $J_{\rho}=[a_{\rho},b_{\rho}]$ of $L[x_{1}',x_{2}']$ by means of geodesics, and $\{a_{\rho}\}\rightarrow a$, $\{b_{\rho}\}\rightarrow b$ when $\rho \rightarrow 0$. In fact, reasoning with $a$, if $a_{\rho}\neq a$ for all $\rho$, let $(c^{\infty}=1,K^{\infty})$ be a point in $\partial[0,1]\times [K^{-},K_{N}(1)]$, which is the limit, up to a subsequence, of $(c^{k},K^{k})\in {\cal C}_{1/k}$. Recall that as the sequence lies in ${\cal C}_{1/k}$ then it must satisfy (ii) (or (iii)) of Lemma \ref{l6} with $a_{\star}^{\infty}=a$, $b_{\star}^{\infty}=b$. If $\int_{\tau_{0}}^{a}(c_{i}^{\infty})^{-1/2}\cdot h_{i}^{\epsilon=-1}[c^{\infty},K^{\infty}]<\infty$ for  $i=2$ then for $k$ big enough each $(c^{k},K^{k})$ corresponds with a fake geodesic and Lemma \ref{l10} can be claimed. Otherwise, the result can be obviously obtained. $\diamondsuit$

\begin{theorem} \label{t6}
Inner Schwarzschild spacetime is geodesically connected even though its corresponding part $(I\times F_{1}, -d\tau^{2}+f_{1}^{2}\cdot g_{1})$ is not.
\end{theorem}

{\it Proof.} The first assertion is a direct consequence of Lemma \ref{lemma'}.

For the second one, let $x_{1}$, $x_{1}'$ be different points in $F_{1}$ and let $l_{1}=d(x_{1},x_{1}')>0$. As Condition (*) does not hold for $c=1$, $i=1$, there exist $\tau_{0}\in I$ such that $\int_{0}^{\tau_{0}}\frac{1}{f_{1}(\tau)}<l_{1}$. From the behaviour of $f_{1}$ towards $\tau_{+}$, for every $\tau_{0}'$ close enough to $0$, the points $(\tau_{0},x_{1})$, $(\tau_{0}',x_{1}')$ cannot be connected by means of a geodesic with $\epsilon=1$. Moreover, from the behaviour of $f_{1}$ towards $0$ and the continuity of the integrals $\int_{a_{\star}}^{\tau_{0}}f_{1}^{-2}(\frac{1}{f_{1}^{2}(\tau)}-D)^{-1/2}$  varying $D$,  for $\tau_{0}'$ close enough to $0$ the points $(\tau_{0},x_{1})$, $(\tau_{0}',x_{1}')$ cannot be connected by means of a geodesic with $\epsilon=-1$ either.  $\diamondsuit$

\begin{remark} \label{r6} {\rm 
It is clear that each two points of inner Schwarzschild spacetime (as well as any multiwarped spacetime for which Lemma \ref{lemma'} is appliable) can be joined by a {\it spacelike} geodesic. From Proposition \ref{p0}, if the second fiber of this spacetime (a round sphere) is replaced by any strongly convex fiber (for example, $\R{}^{2}$) then the resulting multiwarped spacetime is {\it not} geodesically connected. This difference is possible because when the fibers are strongly convex then the geodesic connectedness depends exclusively of the behavior of the warping functions at the extremes of $I$; but if they are just weakly convex, the existence of multiple geodesics in the fibers may maintain the geodesic connectedness under small weakenings to this behavior. 
Thus, even though the conditions for geodesic connectedness in Theorem \ref{t5} are sufficient but not necessary (a necessary condition can be expressed in a reasonably simple way for $n=1$ fiber), the proof of Theorem \ref{t6} and Lemma \ref{lemma'} shows the high accuracy which can be achieved by our technique.  
}
\end{remark}

\end{document}